\documentclass[12pt,reqno]{amsart}

\usepackage{amsmath,amsthm,amscd,amsfonts,amssymb}
\usepackage{eepic}

\usepackage[PostScript=dvips,balance,midshaft,nohug]{diagrams}

\usepackage{tikz, circuitikz}\usetikzlibrary{math}

\swapnumbers  

\newcommand{\QQ}{\mathbb Q}
\renewcommand{\b}{^{\bullet}}

\newtheorem{thm}[subsection]{Theorem}

\newtheorem{lem}[subsection]{Lemma}
\newtheorem{prop}[subsection]{Proposition}
\theoremstyle{definition}

\newcounter{remarkscounter}

\numberwithin{equation}{section}

\makeatletter
\renewcommand\section{\@startsection {section}{1}{\z@}%
                                   {-3.25ex \@plus -1ex \@minus -.2ex}%
                                   {1.8ex \@plus.2ex}%
                                   {\centering\normalfont\normalsize\bfseries}}%
\renewcommand\subsection{\@startsection {subsection}{1}{\z@}%
                                   {-2.0ex \@plus -1ex \@minus -.2ex}%
                                   {-1.8ex \@plus.2ex}%
                                   {\normalfont\normalsize\bfseries}}%
\makeatother
\newcommand{\quash}[1]{}

\newcommand{\ul}[1]{\underline{#1}}
\newcommand{\Hom}{\mathop{\rm Hom}}
\newcommand{\s}{\mathbf A^{\bullet}}
\renewcommand{\t}{\mathbf T^{\bullet}}
\newcommand{\DD}{\mathbb D}
\newcommand{\RR}{\mathbb R}
\newcommand{\CC}{\mathbb C}

\newcommand{\coker}[1]{\mathop{\rm coker}\left(#1\right)}
\newcommand{\B}{\mathcal B}
\newcommand{\rank}{\mathop{\rm rank}}
\newcommand{\uS}{\underline{\bf S}}
 
\renewcommand{\DD}{\mathbb D}
\renewcommand{\SS}{\mathbb S}
\newcommand{\ZZ}{\mathbb Z}
\renewcommand{\c}{I^n} 
\usepackage[utf8]{inputenc}

\begin{document}
\newarrow{Equals} =====
\newarrow{Implies} ===={=>}
\newarrow{Onto} ----{>>}
\newarrow{Into} C--->
\newarrow{Dotsto} ....>

\newcommand{\defeq}{\stackrel{\rm def}{=}}

\title{Topological aspects of Boolean functions}
\author{\rm Anders Bj{\"{o}}rner}
\address{Department of Mathematics, KTH Stockholm}
\author{Mark Goresky}
\address{School of Mathematics, IAS, Princeton}
\author{Robert MacPherson}\address{School of Mathematics, IAS, Princeton}

\dedicatory{To Claudio Procesi, with friendship and admiration, on the occasion of his eightieth birthday}

\begin{abstract}
We discuss ways in which tools from topology can be used to
derive lower bounds for the circuit complexity of Boolean functions.
\end{abstract}
\keywords{Boolean function, circuit complexity, order complex, Betti number, simplicial sheaf}
\subjclass{94D10, 94C11, 06A07, 13F55}
\maketitle

\setcounter{tocdepth}{1} %
\tableofcontents

\section{Introduction}  By a {\em Boolean function} (of $n$ variables) we mean a function $f:\{0,1\}^n\to \{0,1\}$.  
Various notions of the computational complexity of a Boolean function $f$ exist,
 but all are asymptotically equivalent to the
 {\em circuit complexity} $CC(f)$,  (or {\em formula size}, see \S \ref{sec-circuits}), which is notoriously difficult to compute. 
 It is known (see, e.g., \cite{Circuits}) that random
Boolean functions tend to have exponential circuit complexity.  

It is an interesting goal to construct
for any Boolean function $f$ a topological space $X_f$,
 with the property that some standard topological invariants
 applied to $X_f$ would give estimates 
for the circuit complexity of $f$. 
Such an approach
 has been tried with some success in
the case of some non-Boolean functions.
See e.g. \cite{Bjorner5}, where the topological model  is
arrangements of
linear subspaces and the topological measure is
the sum of Betti numbers.

In this paper we take the first steps towards the goal of exploring 
the potential use of ``topological complexity''  in the setting of Boolean functions $f$.
We introduce two constructions of a topological space associated to
$f$:  the {\em red corner model} $\Sigma(f)$, obtained from the order complex of the partially orderd set defined by $f^{-1}(1)$, 
 and the {\em subspaces model} $Z(f)$, a collection of coordinate subspaces and their complements in $\CC^n$, defined by $f$. 
 For these  models we define,  in equations (\ref{eqn-beta}) and (\ref{eqn-gamma}), certain
 integers $\beta(f)$ and $\gamma(f)$, essentially the sum of the Betti numbers of the corresponding
 topological space.  We prove several basic computational facts about these invariants.

Using the standard rules of logic, Boolean
functions may be combined using $\wedge$ (``AND"), $\vee$ (``OR") and $\sim$ (``NOT"), leading to computational issues.
The permutation group $\mathfrak{S}_n$ on $n$ letters,  acts on the set of Boolean functions by setting 
\[ f^{\sigma}(x_1,x_2,\cdots,x_n) = f(x_{\sigma(1)}, x_{\sigma(2)},\cdots, x_{\sigma(n)})\]
for any $\sigma \in \mathfrak{S}_n$.  
If $f: B_n \to \{0,1\}$ and $g: B_m \to \{0,1\}$ are Boolean functions, we use the identification
$B_n \times B_m \to B_{n+m}$ to define Boolean functions $f\wedge g$ and $f\vee g$.
\begin{thm}\label{thm-intro}
Suppose $f:\{0,1\}^n \to \{0,1\}$ is a Boolean function. Define $\beta(f)$ as in (\ref{eqn-beta}) and $\gamma(f)$
as in (\ref{eqn-gamma}).   Then 
\begin{enumerate}
\item\label{prop-symmetric} $\beta(f^\sigma) = \beta(f)$ and $\gamma (f^\sigma) = \gamma(f)$, for any $\sigma$ in the symmetric group $\mathfrak S_n$,
\item\label{prop-duality}
$\beta(\sim f) = \beta(f)$.  If $f$ is not identically $0$ or $1$ then 
$\gamma(\sim f) = \gamma(f)$.
\item \label{prop-sum}
$\beta(f \vee g) = \beta(f \wedge g) = \beta(f)\beta(g)$.
If neither $f$ nor $g$ is identically $0$ or $1$, then $\gamma(f \vee g) = \gamma(f \wedge g) = \gamma(f)\gamma(g)$.
\item \label{prop-zero}
Suppose $f$ does not depend on one of the variables:  $f(x_1,x_2,\cdots,x_n) = g(x_2,\cdots,x_n)$. 
Then $\beta(f) = 0$ and $\gamma(f) = \gamma(g)$.
\item\label{prop-tree} Suppose the Boolean function $f$ may be computed by a Boolean circuit (see \S \ref{sec-circuits}) 
that is a tree.  If $f$ is not identically $0$ or $1$, then $\beta(f) = 1$ and $\gamma(f) \le 2^n$.
\end{enumerate}
\end{thm}
 This is proved in \S \ref{sec-RedCornerProof} and \S \ref{sec-proof-subspaces}, respectively.
In \S \ref{sec-speculation} we suggest additional techniques to associate a space with a Boolean function.

\subsection*{Acknowledgements} The second and third authors are grateful for the hospitality shown by Institut Mittag-Leffler and KTH, Stockholm,
during their visits in 2008 and 2013. The authors are grateful to Richard Ehrenborg and an anonymous referee for many useful comments and suggestions.
This paper was typeset by \AmS-\LaTeX  using Paul Taylor's package diagrams.sty 

\section{Boolean functions and  circuits}\label{sec-circuits}
\subsection{Boolean circuits}
Identify $0$ with ``False" and $1$ with ``True" and denote $\vee=$ ``OR" and $\wedge=$ ``AND" and $\sim$ = ``NOT".
By a {\em Boolean circuit} with $n$ inputs we shall mean a circuit consisting
of AND and OR gates having a fan-in of two, a fan-out of one, and
arbitrarily many NOT gates; the whole circuit having $n$ inputs and one
output.   (By deMorgan's rule, all
NOT gates may be moved to the bottom of the circuit.) 
The output of such a circuit $\Omega$ defines a Boolean function  $f_{\Omega}:
\left\{0,1\right\}^n \to \left\{0,1\right\}$. We say that the circuit $\Omega$
computes the function $f_{\Omega}$. 
 
\begin{figure}[h!]
\begin{circuitikz} \draw
(2,0) node[and port, rotate=90] (gate1) {$>$}
(0,0) node[or port,rotate=90] (gate2) {$<$}
(1,2) node[and port, rotate=90] (gate3) {$>$}
(gate1.out) to[isourcesin, sources/scale=0.5]  (gate3.in 2)
(gate2.out) -- (gate3.in 1);
\draw (0,-2) node[mixer,scale=0.3,label=below:$x_1$](xin) {};
\draw (xin) -- (gate1.in 1); \draw (xin) -- (gate2.in 1);
\draw (2,-2) node[mixer,scale=0.3,label=below:$x_2$](yin){};
 \draw(yin)--(gate1.in 2);\draw(yin) -- (gate2.in 2);
\draw (gate3.out) -- (1,2.5);
\draw [-latex] (1,2.5) node[label=above:$(x_1+x_2)$ mod 2 ]{};
\end{circuitikz}
\caption{A circuit to calculate ``parity": the sum modulo $2$}
\end{figure}
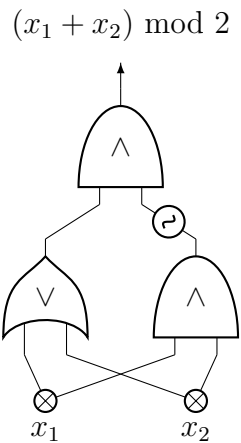
\medskip\noindent
 Denote by $|\Omega|$ the number of AND/OR
gates in the circuit.  (We do not count the NOT gates.)  The {\em circuit complexity}
of a Boolean function $f$ is 
\[CC(f) \defeq \min \{ |\Omega| |\ \Omega \text{ computes } f \}. \]

A theorem of C. Shannon (see \cite{Circuits} \S 1.5) states (roughly) for large $n$ that 
most Boolean functions with $n$ inputs require circuits of exponential size, $2^n/n$.

\subsection{Product of valences}
\label{subsec-conjecture} 
For each $1 \le i \le n$, let $a_i$ denote the {\em valence} or fan-out of
the input node $x_i$, so that $\sum_{i=1}^n a_i$ is the number of leaves
on the tree that is obtained by cutting off all the inputs.  
Define the number $PV(\Omega) = \prod_{i=1}^n a_i,$ the product of the
valences for each of the variables, and define $PV(f)$ to be the minimum of $PV(\Omega)$ over
all circuits $\Omega$ that compute $f$.  A Boolean circuit $\Omega$ is a tree if and only if $PV(\Omega)=1$.

The number of gates is one less than the number of  inputs.
Hence, the arithmetic-geometric mean inequality implies that 
\begin{equation}\label{eqn-geometric-mean}CC(f)+1 \ge n(PV(f))^{1/n},\end{equation}
showing that the product-of-valences number  $PV(f)$ has the potential
to lower-bound circuit complexity.  
\subsection*{Conjecture}  Suppose the Boolean circuit $\Omega$ computes $f$.  Then
\begin{equation}\label{eqn-conjecture} PV(\Omega) \ge \gamma(f).\end{equation}
See §8 for motivation for this conjecture.

Computer experiments suggest that for each $n$, the parity function maximizes $\beta(f)$ among all Boolean functions $f$ of $n$ variables, and similarly for $\gamma(f)$.
In particular, even if the conjecture is true, 
these models fail to detect exponential complexity.

\subsection{Boolean lattices}\label{subsec-definitions}
The {\em Boolean lattice} $\B_n$  is the set
 of subsets $I\subset[n]=\{1,2,\cdots,n\}$, partially ordered by inclusion and
ranked by the cardinality $|I|$.  We identify $\B_n \cong \{0,1\}^n$,
associating a subset $I \subset [n]$ with the Boolean vector $x_I = (x_1,x_2,\cdots,x_n)$ where $x_i = 1$ if and only if  $i \in I$.
If $x = (x_1,x_2,\cdots, x_n) \in \B_n$ then its rank is $|x| = \sum_{i=1}^nx_i$, and 
 $x \le y$ if and only if  $x_i \le y_i$ for all $1 \le i \le n$.
 The unique minimal element $\ul{\mathbf 0}\in \B_n$ corresponds to the empty set  and the maximal element $\ul{\mathbf 1}$ corresponds to
 $[n] = \{ 1, 2, \cdots, n\}$.
  The {\em punctured Boolean lattice}
 $\B'_n$ is the lattice $\B_n$ with the minimal element $\ul{\mathbf 0}$ removed.
 If $m,n \ge 1$  the Boolean lattice $B_{m+n}$ may be identified with
the partially ordered set $B_m \times B_n$ (where $(v,w) \le (v',w') \iff v \le v' \text{ and }
w \le w'$). 

\subsection{Sum of Betti numbers}  If $P$ is a partially ordered set (poset) we denote its order complex by $\Sigma(P)$.
It is the abstract simplicial complex with one $k$-simplex for each chain $v_0 < v_1 < \cdots < v_k$ of elements of $P$
and the obvious face relations.
If $\Sigma$ is a simplicial complex we denote by $|\Sigma|$ its geometrical realization.  In this paper
every simplicial complex has a canonical geometrical realization.  If $X$ is a topological space denote by
\[b^*(X) \defeq \sum_{i \ge 0} \text{rank } H_i(X;\QQ)\]
the sum of its Betti numbers.  Denote by $\tilde b^*(X)$ the sum of Betti numbers of the reduced
homology of $X$ with the convention that $\tilde{H}_{-1}(\phi; \QQ) = \QQ$.
 If  $\mathcal M$ is a sheaf of rational vector spaces on $X$,
 denote by $b^*(X,\mathcal M)$ the
sum of the ranks of its cohomology groups. 

\subsection{Order complex}\label{subsec-punctured-order-complex}
If $f:\B_n \to \{0, 1\}$  is a Boolean function,  let $\B(f) = f^{-1}(1)$ denote the poset of elements $v$
such that $f(v) = 1$ with its induced partial order. 
The order complex $\Sigma(f) =\Sigma(\B(f))$ is the full subcomplex of
$\Sigma(\B_n)$ with one $k$-simplex for every chain $x_0<x_1\cdots <x_k$ of elements such that each $f(x_j) = 1$. 
Denote its canonical geometric realization by $|\Sigma(f)|  \subset [0,1]^n$,
where  $[0,1]$ is the real unit interval.

 If $f:\B'_n \to \{0,1\}$ is a {\em punctured} Boolean function (or equivalently, if $f(\ul{\mathbf 0}) = 0$)
 then using (\ref{eqn-simplicial-isomorphisms}) 
 its order complex  may be considered as a subcomplex of $\Delta'$,  hence
 $\Sigma(f)\subset \Sigma(\B'_n)\cong \Delta'$.

\subsection{Geometry of the cube}\label{subsec-cube}
The  order complex $\Sigma(\B_n)$ of the Boolean lattice has a canonical geometrical realization  
\[|\Sigma(\B_n)| \cong [0,1]^n=\c\] 
as a triangulation of the $n$ dimensional cube $\c$. Namely, 
let $e_1,e_2, \cdots , e_n$ denote the standard basis of $\RR^n$,
and for each subset $J \subset [n] = \{ 1, 2, \cdots, n\}$ let 
$$\widehat{e}_J \defeq \sum_{j \in J}e_j \in [0,1]^n$$ and set
 $\widehat{e}_{\phi} \defeq 0$.  
 The vectors $\{\widehat{e}_J\}_{J\subset [n]}$  form the set of
vertices of the $n$-dimensional cube $[0,1]^n$ and correspond to the elements of the Boolean lattice.  
The resulting triangulation of the cube $\c$ has one $k$-simplex
\[ \sigma = \left< v_{J_0}, v_{J_1}, \cdots, v_{J_k}\right>\]
for each chain of subsets $J_0 \subset J_1 \subset \cdots \subset J_k \subset [n]$.
There are $2n$ facets (codimension one faces) of the cube:  for each choice $1 \le i \le n$ 
  \[ \partial_i^- \c = \left\{ (x_1,\cdots,x_n)\in \c \left| \  x_i = 0 \right. \right\} \]
  \[ \partial_i^+ \c = \left\{ (x_1,\cdots,x_n)\in \c \left| \ x_i = 1 \right. \right\} \] 
Each facet $\partial_i^{-}\c$ contains the origin $\ul{\mathbf 0} = (0,0, \cdots,
 0)$, and each facet $\partial_i^{+}\c$
contains the vertex  $\ul{\mathbf 1} = (1, 1, \cdots, 1)$. 
Let \[\partial^+\c = \bigcup_{i \in [n]} \partial^+_i\c.\]  This  space is the union of all facets in $\c$ that do not contain
the origin.\medskip

\newcommand*{\haoffset}{.2}
\newcommand*{\hboffset}{.4}
\newcommand*{\hcoffset}{.6}

\begin{figure}[h!]
\begin{tikzpicture}[scale = 2.5]

\def\del{(\hboffset - \haoffset)/2};

\colorlet{mycolor}{red!20};
\node (O) at (0,0) {$\phi$};
\node (A1) at (\haoffset-1, 1){\{1\}};
\node(A2) at  (\haoffset,1){\{2\}};
\node (A3) at (\haoffset +1,1){\{3\}};

\node (B1) at (\hboffset -1, 2) {\{12\}};
\node (B2) at (\hboffset, 2) {\{13\}};
\node (B3) at (\hboffset+1,2) {\{23\}};\node (C1) at (\hcoffset,3) {\{123\}};

\draw [very thick](O) -- (A1); \draw[very thick] (O) -- (A2); \draw [very thick](O) -- (A3);

\draw [very thick] (A1) -- (B1);  \draw [very thick](A3) -- (B3);
\draw[very thick, dotted] (A1) -- ({\haoffset + \del -.55}, 1.45);
\draw[very thick, dotted] ({\haoffset +\del -.45}, 1.55) -- (B2);
\draw[very thick, dotted] (A3) -- ({\haoffset + \del +.55}, 1.45);
\draw[very thick, dotted] ({\haoffset +\del +.45}, 1.55) -- (B2);

 \draw [very thick] (A2) -- (B1);
\draw[very thick] (A2) -- (B3);
 
  \draw [very thick] (B1) -- (C1); \draw [very thick, dotted] (B2) -- (C1); \draw [very thick](B3) -- (C1);

\shadedraw[ultra nearly transparent, fill=red!10!black!10] (\haoffset +1,1) -- (\hboffset,2) -- (\hcoffset,3) -- (\hboffset+1,2) -- (\haoffset+1,1); 
\shadedraw[ultra nearly transparent, fill=green!1!black!1] (\haoffset-1,1) -- (\hboffset -1,2) -- (\hcoffset,3) --(\hboffset,2) -- (\haoffset-1,1);
\shadedraw[ultra nearly transparent, fill = green!10!black!10!] (0,0) -- (\haoffset -1,1) -- (\hboffset,2) -- (\haoffset +1, 1)-- (0,0);

\node (E2) at (3,1) {\{2\}};
\node (E1) at (2,3) {\{1\}};
\node (E3) at (4,3) {\{3\}};
\node (E12) at (2.5,2) {\{12\}};
\node (E23) at (3.5,2) { \{23\} };
\node (E13) at (3,3) { \{13\} };
\node (E123) at (3, 2.5) {\{123\}};

\draw [very thick] (E2) -- (E12);
\draw [very thick] (E12) -- (E1);
\draw [very thick] (E2) -- (E23);
\draw [very thick] (E23) -- (E3);
\draw [very thick, dotted] (E1) -- (E13);  \draw [very thick, dotted] (E13) -- (E3);
\draw [red,dotted] (E123) -- (E1); \draw [red, dotted] (E123) -- (E2); \draw [red, dotted] (E123) -- (E3);
\draw [very thick] (E123) -- (E12); \draw [very thick, dotted] (E123) -- (E13); \draw [very thick] (E123) -- (E23);
\end{tikzpicture}
\caption{The cube $\Sigma(\B_3)$ and the simplex $\Delta'$}
\vspace{8mm}
\end{figure}
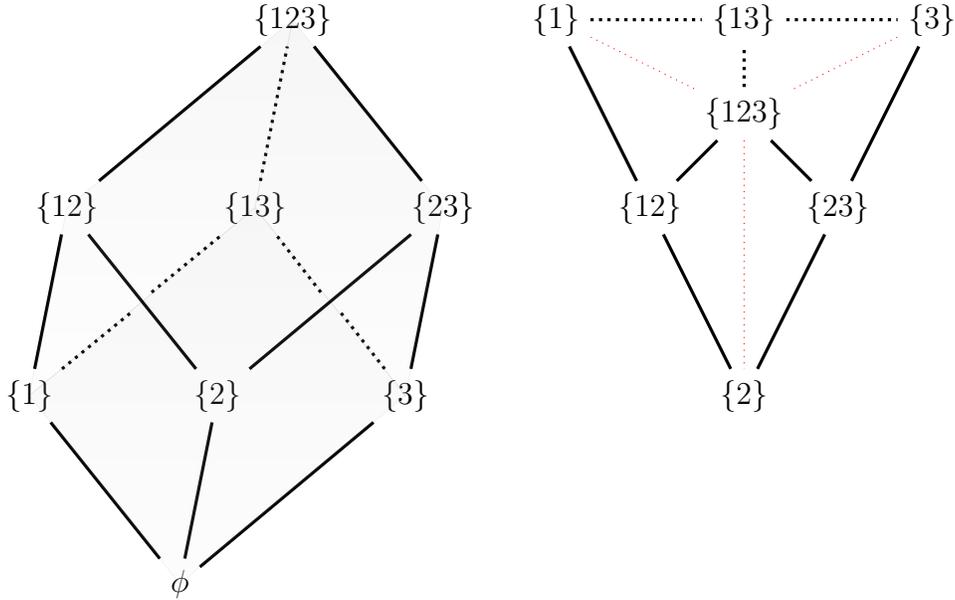

\subsection{Geometry of the simplex}\label{subsec-defs2}
 The $n-1$ dimensional simplex
\[\Delta = \Delta^{n-1} \defeq  \left\{ (x_1, x_2, \cdots, x_n) \in  \RR^n
\left| 0 \le x_i \le 1, \ \sum_{i=1}^n x_i = 1 \right. \right\} \]
is the convex hull of its vertices, the
standard basis vectors $e_1, e_2, \cdots, e_n$.  Each nonempty subset $J \subset [n]$ determines
a face $e_J\defeq \text{convex span} \left\{ e_j | \ j\in J\right\}$.  Each $j \in [n]$ determines a {\em facet}
\[\partial_j \Delta \defeq 
\Delta \cap \{ x \in \RR^n \left| \ x_j = 0 \right. \}= e_{[n]-\{j\}}.\]

Let $\Delta'$ denote the barycentric subdivision of $\Delta$. 
 It is the order complex
 \[ \Delta' = \Sigma (\B'_n)\]
of the poset $\Delta =\B'_n$ of non empty faces of $\Delta$.   It has one vertex $\widehat{e}_J$ for
each (non empty) face $e_J$ of $\Delta$. To be explicit, the correspondence
\begin{diagram}[size=2em]
\text{nonempty subsets of }[n]\ &\longleftrightarrow &
\text{\ faces of\ }\Delta\ & \longleftrightarrow &
\text{\ vertices of \ } \Delta' & \longleftrightarrow & 
\text{\ nonzero Boolean vectors}\\
J &\longleftrightarrow & e_J & \longleftrightarrow & \widehat{e}_J & \longleftrightarrow& (x_1,x_2,\cdots,x_n)
\end{diagram}
(where $x_j = 1 \iff j \in J$) determines simplicial isomorphisms
\begin{equation}\label{eqn-simplicial-isomorphisms}
 \left|\Delta'\right| \longleftarrow \left| \Sigma(\B'_n) \right| \longrightarrow \partial^+\c.\end{equation}
 
 The simplicial complex described in various ways by 
(\ref{eqn-simplicial-isomorphisms}), 
is topologically an $(n-1)$-dimensional ball. More precisely, 
 it is the cone over an $(n-2)$-dimensional sphere.

In \S \ref{subsec-partial_J} we shall use the following observation:
 If $\tau: \partial^+\c \to |\Delta'|$ denotes the composed isomorphism in (\ref{eqn-simplicial-isomorphisms})
 then for any $j \in [n]$, 
 \begin{equation}\label{eqn-facet-intersection}
 \tau\left(\partial_j^- \c \cap \partial^+\c\right) = \partial_j\Delta.\end{equation} 

\section{The red corner model}
For each vertex $v\in \{0,1\}^n$ let $R_v$ be the union of those facets (codimension one faces)
$F\subset \c$ such that $v \in F.$  We sometimes refer to $R_v$ as
the {\em red corner} at $v.$ The remaining facets of $\c$ are those
in $R_{\sim v}$, the red corner of the opposite vertex.
The {\em total Betti number} $\beta(f)$ in the red corner model is 
\begin{equation}\label{eqn-beta}
\beta(f)  \defeq 
\sum_{v\in \B_n} b^*(|\Sigma(f)|, R_v\cap|\Sigma(f)|).
\end{equation}
It is a sum of $2^n$ terms, each of which is the sum of all the Betti
numbers of the relative homology of the pair $(|\Sigma(f)|, R_v\cap|\Sigma(f)|).$

\subsection{Proof of Theorem \ref{thm-intro} for the red corner model}\label{sec-RedCornerProof}
Part (\ref{prop-symmetric}) is clear from the definition.  Parts (\ref{prop-duality}) and
(\ref{prop-sum}) will be verified in
\S \ref{subsec-duality} and \S \ref{subsec-sum}.  For part (\ref{prop-zero}), if the function $f$ does
not depend on a variable, say $x_n$ then $f(x_1,\cdots,x_n) = f'(x_1,\cdots,x_{n-1})$ for some
Boolean function $f'$ of $n-1$ variables, so $\Sigma(f)$ is the Cartesian product of $\Sigma(f')$ and
the unit interval.  This implies that for every vertex $v \in \c$ and for every simplex $\sigma$ of
$\Sigma(f')$ either $\sigma\times\{0\} \subset R_v$ or else $\sigma \times {1} \subset R_v$.  Therefore
the relative homology $H_*(|\Sigma(f)|, R_v\cap|\Sigma(f)|)$ vanishes.  Part (\ref{prop-tree}) follows from part
(\ref{prop-sum}) by induction: at the top gate  each branch involves different
input variables. 

\subsection{Sheaves and the red corner model}\label{sec-duality}
 Throughout this paper the word ``sheaf" refers
to a sheaf of finite dimensional $\QQ$-vector spaces.

On the circle $T=S^1$ let  $m$ denote the M\"obius local system.  It is the rank one
local coefficient system over $\QQ$ such that the monodromy around the circle is multiplication by $-1$.
One checks that $
H^i(S^1, m) = 0 \text{ for all } i \ge 0 $
and that multiplication $ m \otimes  m \to \underline{\QQ}$ induces a {\em self duality} isomorphism  $  m \cong \Hom( m, \underline{\QQ}),
$
where $\underline{\QQ}$ denotes the trivial (rank one) local system on $S^1$.  

Let $r_0$ be the
sheaf on the unit interval $[0,1]$ that is the constant sheaf $\underline{\QQ}_{(0,1]}$ on 
$(0,1]$ and is zero on the point $\{0\}$.  Let $r_1$ be the sheaf that is 
$\underline{\QQ}_{[0,1)}$ on $[0,1)$ and is zero on $\{1\}$.

\begin{lem}\label{lem-circle}
Let $\pi:S^1 \to [0,1]$ be the mapping $\pi(e^{i\theta}) = {\scriptstyle\frac{1}{2}}(1+cos(\theta))$.
Then there is an isomorphism of sheaves, $\pi_*(m) \cong r_0 \oplus r_1$. 
 \end{lem}
 
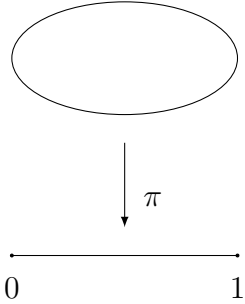
\begin{figure}[h!]
\begin{tikzpicture}[scale=0.75]
\draw (5,5) ellipse (2 cm and 1 cm);
\draw[-latex] (5,3.5) -- (5,2) {}; \draw(5.5,2.5) node{$\pi$};
\draw (3,1.5) node[label=below:$0$](zero){};
\draw(7,1.5) node[label=below:$1$](one){};
\draw (3,1.5) -- (7,1.5);\draw[fill](3,1.5) circle [radius = .025];
\draw[fill](7,1.5) circle[radius=.025];
\end{tikzpicture}
\caption{The mapping $\pi$}
\end{figure}
\begin{proof}
A simplicial sheaf $A$ on the interval $[0,1]$ is an arrangement of $\QQ$-vector spaces and morphisms,
\begin{diagram}[height=1em] A_{\{0\}}  &\rTo& A_{(0,1)} & \lTo & A_{\{1\}}. \end{diagram}  An isomorphism of simplicial
sheaves, $\Phi: \pi_*(m) \to r_0 \oplus r_1$ is defined by the following  commutative diagram,
where $\Phi(a,b) = \left(\frac{a-b}{2}, \frac{a+b}{2}\right)$:
\begin{diagram}[height=1em]
\pi_*(m): \qquad & \QQ & \rTo & \QQ \oplus \QQ & \lTo & \QQ & \\
         & x & \rTo & (x,-x) &&&\\
          & & & (y,y) & \lTo & y&\\   &&& {} &&&\\ \\
          &&&\dTo_{\Phi} && & \\ &&& {} &&& \\  
  r_0: \qquad & \QQ & \rTo^{id} & \QQ & \lTo & 0 &\\
  &&&\oplus &&&\\
  r_1: \qquad & 0 &\rTo & \QQ & \lTo_{id} & \QQ &\qedhere
\end{diagram}
 \end{proof}

\subsection{The torus.}
On the torus $T^n=(S^1)^n$ let $\mathcal M = m\boxtimes m\boxtimes\cdots\boxtimes m$ be the tensor
product of the M\"obius local systems from the factors.  It follows from the K\"unneth theorem
that  $H^i(T^n,\mathcal M)=0$ for all $i \ge 0$, and that $\mathcal M \cong {\rm Hom}(\mathcal M, \underline\QQ)$ 
is self dual.  Let $\pi:T^n \to [0,1]^n$ be the
product of the projections.  It follows that there is an isomorphism of sheaves,
\[ \pi_*(\mathcal M) \cong \underset{v\in \B_n}{\bigoplus} \mathcal M_v\]
where, for each vertex $v$ of the cube $\c$, the sheaf $\mathcal M_v$ is zero on $R_v$
(the red corner corresponding to $v$) and is the constant sheaf $\underline{\QQ}$ on the complement, 
$\c - R_v$. So the cohomology of $\mathcal M_v$ is the relative cohomology, that is, for any
closed subset $W \subset I^n$ we have: $H^i(W, \mathcal M_v) = H^i(W, W\cap R_v)$.

 Since $\pi$ is a finite mapping, we conclude that for any Boolean function
$f$, and for any $i \ge 0$ there is an isomorphism (see Remark \ref{subsec-rmk} below),
\begin{align*}
 H^i(\pi^{-1}(|\Sigma(f)|), \mathcal M)
&\cong H^i(|\Sigma(f)|, {\rm R}\pi_*\mathcal M) \cong H^i(|\Sigma(f)|, \pi_*(\mathcal M))\\
&\cong \bigoplus_{v \in \B_n}H^i
(|\Sigma(f)|, \mathcal M_v) \cong
\bigoplus_{v\in \B_n} H^i(|\Sigma(f)|,R_v \cap |\Sigma(f)|).
\end{align*}
Hence $\beta(f) = \sum_{i \ge 0} \text{rank} H^i(\pi^{-1}(|\Sigma(f)|), \mathcal M)$.

\subsection{Remark.}\label{subsec-rmk}
Suppose $h:X\to Y$ is a continuous map between locally compact Hausdorff spaces and suppose that
$S$ is a sheaf on $X$ whose stalk cohomology is finite dimensional at each point.  Let
$h_*(S)$ be its pushforward to $Y$.  Fix $i \ge 0$.  It is not always
true that $H^i(X,S) \cong H^i(Y, h_*(S))$.  However, there is a canonical
isomorphism $H^i(X,S) \cong H^i(Y, {\rm R}h_*(S))$ where ${\rm R}h_*(S)$ is the push-forward of the sheaf $S$,
in the sense of the derived category of constructible sheaves.  In general, ${\rm R}h_*(S)$
will be a complex of sheaves, rather than a single sheaf.  However, in our case 
the mapping $\pi$ is finite and it is even simplicial, so there is an isomorphism in the 
derived category between $\pi_*(\mathcal M)$ and ${\rm R}\pi_*(\mathcal M)$.  More
classically, one would say that the Leray-Serre spectral sequence for the map $\pi$
collapses because the $E^2$ page lives on a single line.

We will need the following standard facts.
\begin{lem}\label{lem-P1}
Let $X$ be a simplicial complex with vertex set $V$.  Let $\Sigma_0 \subset X$ be a full subcomplex,
spanned by a subset of the vertices $V_0\subset B$.  Let $\Sigma_1 \subset X$ be the {\em supplementary}
subcomplex, that is, the full subcomplex spanned by the set of remaining vertices $V_1 = V - V_0$.  Then
there is a simplicial map $h:X \to \Delta^1$ {\rm (}so $h:|X| \to [0,1]${\rm )} with $h^{-1}(t) = \Sigma_t$ $(t = 0,1)$ and there are
natural deformation retractions, $|X| - |\Sigma_0| \to |\Sigma_1|$ and $|X| -|\Sigma_1| \to |\Sigma_0|$.  \qed
\end{lem}
\begin{lem}\label{lem-P2}
Let $X$ be a smooth compact $n$-dimensional manifold and let $\mathcal M$ be a local 
coefficient system of $\QQ$-vector spaces on
$X$.  Let $\mathcal M^*$ be the dual local coefficient system.  Let $\Sigma \subset X$ be a
closed subset that is a subcomplex with respect to some smooth triangulation of $X$.
Then the pairing $\mathcal M \times
\mathcal M^* \to \underline{\QQ}$ induces the Poincar\'e duality isomorphism
\[ H^j(X, X-\Sigma, \mathcal M) \cong H_{n-j}(\Sigma, \mathcal M^*)
\cong \left(H^{n-j}(\Sigma, \mathcal M^*)\right)^*.  \tag*{\qed}\]
\end{lem}
\subsection{Proof of Theorem \ref{thm-intro}(part \ref{prop-duality})}  \label{subsec-duality}

Let $g =\ \sim\!\!f$ be the negation of $f$.  Then $\Sigma(g)$ is supplementary to $\Sigma(f)$ in
the simplicial complex $\c$ so $\pi^{-1}(\Sigma(f))$ is supplementary to $\pi^{-1}(\Sigma(g))$ in the torus $T$
and there is a deformation retraction
\[ T^n - \pi^{-1}(|\Sigma(f)|) \to |\Sigma(g)|.\]
Using Lemma \ref{lem-P1}, Lemma \ref{lem-P2} and the fact that the sheaf $\mathcal M$ is acyclic and self-dual on $T^n$, 
it follows, for each $i \ge 0$ that
\begin{align*}
 H^i(\pi^{-1}(|\Sigma(g)|),\mathcal M) &\cong H^i((T^n - \pi^{-1}|\Sigma(f)|,\mathcal M)\\
&\cong H^{i-1}((T^n, T^n - \pi^{-1}|\Sigma(f)|, \mathcal M)\\
&\cong H_{n-i+1}(\pi^{-1}|\Sigma(f)|, \mathcal M^*)\\
&\cong H_{n-i+1}(\pi^{-1}|\Sigma(f)|, \mathcal M)\\
&\cong H^{n-i+1}(\pi^{-1}|\Sigma(f)|, \mathcal M)^*. \end{align*}
where ``$*$" denotes the (rational vector space) dual.  Consequently
\[ \beta(g)=\sum_{i \ge 0} \text{rank}H^i(|\Sigma(g)|,\mathcal M) = 
\sum_{j \ge 0} \text{rank}H^j(|\Sigma(f)|, \mathcal M)=\beta(f).  \qedhere\]
This completes the proof of part \ref{prop-duality} of Theorem \ref{thm-intro}.

\subsection{Proof of Theorem \ref{thm-intro}(part \ref{prop-sum})}\label{subsec-sum}
Let $f:B_n \to \{0,1\}$ and $g:B_m \to \{0,1\}$ be Boolean functions.  Using the identification
$B_n \times B_m \to B_{n+m}$ we have that
\begin{align*} 
\Sigma(f \wedge g) &= \Sigma(f) \times \Sigma(g)\\
\Sigma(f \vee g) &= \left(\Sigma(f) \times I^m\right) \cup 
\left(\c \times \Sigma(g)\right).
\end{align*}
Lift these statements to the torus 
\[T^n \times T^m \cong T^{n+m} \underset{\pi_{n+m}}{\longrightarrow} I^{n+m}
\] where we have an
isomorphism of M\"obius local systems $\mathcal M_n \boxtimes \mathcal M_m \cong \mathcal M_{n+m}$. 
Using the K\"unneth formula gives
\[ H^*(\pi_{n+m}^{-1}(|\Sigma(f\wedge g)|,\mathcal M_{n+m}) 
\cong H^*(\pi_n^{-1}(|\Sigma(f)|,\mathcal M_n) \otimes H^*(\pi_m^{-1}
|\Sigma(g)|,\mathcal M_m)\]
and
\begin{equation}\label{eqn-vanishing}
 H^*(\pi_{n+m}^{-1}(|\Sigma(f)| \times I^m),\mathcal M_{n+m}) = 
 H^*(\pi_{m+n}^{-1}(\c \times |\Sigma(g)|,\mathcal M_{n+m}) = 0.\end{equation}
The first equation says that $\beta(f \wedge g) = \beta(f)\beta(g)$.  For $f \vee g$ we use the Mayer
Vietoris theorem for the two sets $|\Sigma(f)| \times I^m$ and $\c \times |\Sigma(g)|$ whose
intersection is $|\Sigma(f \wedge g)|= |\Sigma(f)| \times |\Sigma(g)|$.  From the vanishing
condition (\ref{eqn-vanishing}) it follows that $\beta(f \vee g) = \beta(f\wedge g)= \beta(f) \beta(g)$.
(Alternatively, one may use Alexander duality to reduce the case of $f \vee g$ to that of $f \wedge g$.)  \qed

\subsection{Computations with red corners}

With notation as in the preceding section, observe that
\begin{align*}
R_v \cap \Sigma(f) &= \{(x_0<\cdots<x_k)\in \Sigma(f)|\ x_0 \wedge v \ne \ul{\mathbf 0} \text{ OR }  
x_k \vee v \ne \ul{\mathbf 1}  \}\\
&=\Gamma_v^+(f) \cup \Gamma_v^-(f)\end{align*}
where
\begin{align*}
\Gamma_v^+(f) = \{ (x_0<x_1<\cdots<x_k)\in \Sigma(f)|\ x_0 \wedge v \ne \ul{\bf 0} \} \\
\Gamma_v^-(f) = \{ (x_0<x_1<\cdots<x_k) \in \Sigma(f)|\ x_k \vee v \ne \ul{\bf 1} \} \end{align*}
The complex $R_v \cap\Sigma(f)$ is not, in general, an order
complex:  There may be $y_1<\cdots <y_k$ and $y_k < \cdots < y_m$ both belonging to $R_v \cap \Sigma(f)$ but
$y_1 <\cdots < y_k \cdots < \cdots y_m$ not belonging. However, the following is an order complex, which facilitates
computation with known methods:
\begin{align*}
\Gamma_v(f) &= \Gamma_v^+(f) \cap \Gamma_v^-(f)\\
&= \{(x_0<\cdots<x_k)\in \Sigma(f)|\ x_0 \wedge v \ne \ul{\mathbf 0} \text{ AND } x_k \vee v \ne \ul{\mathbf 1}\} \\
& = \Sigma \left (\{ x \in \B_n|\ f(x) =1 \text{ AND }   x \nless \sim v \text{ AND  } x \ngtr \sim v \}
\right) \end{align*}
where $\sim v$ denotes the complementary 0/1 vector in $\B_n$. 

\begin{prop}\label{prop-Gammav}
  Let $f:\B_n \to \{0,1\}$ be a Boolean function with\footnote{This is not a serious restriction.
If $f(\ul(\mathbf 0)) = f(\ul(\mathbf 1)) = 1$ use the opposite Boolean function.  Otherwise
changing either of these values will have minimal effect on the complexity of $f$.} $f(\ul{\mathbf 0}) = f(\ul{\mathbf 1}) = 1$. 
 Then  
\[ \beta(f) = \sum_{v\in\B_n} \tilde{b}^*(\Gamma_v(f)).\qedhere\]
\end{prop} 
\noindent
 Here $\tilde b^*(X)$ denotes the sum of the Betti numbers of the {\em reduced} 
 homology of $X$, with the convention that $H_{-1}(\phi; \ZZ) = \ZZ$.

\begin{proof}
 Since $\Sigma(f)$ is contractible, 
\[H_i(\Sigma(f), R_v \cap \Sigma(f)) \cong \tilde{H}_{i-1}(R_v \cap \Sigma(f)) = \tilde{H}_{i-1}(\Gamma_v^+(f) \cup \Gamma_v^-(f)).\]
Both $\Gamma_v^+(f)$ and $\Gamma_v^-(f)$ are contractible so Mayer-Vietoris for reduced homology  gives
\[ \tilde{b}^*(R_v \cap \Sigma(f)) = \tilde{b}^*(\Gamma_v(f)).\qedhere\]\end{proof}

\section{The subspaces model}\label{sec-subspaces}
\subsection{Stratification}
The Boolean lattice $\B_n$ may be used to define a stratification of $\RR^n$ consisting of coordinate subspaces and their complements. 
To each Boolean vector $x=(x_1,\cdots,x_n)\in \B_n$  we associate the following ``stratum" in $\RR^n$:
\[ Z_x \defeq \left\{ z \in \RR^n\left|\
z_i = 0 \iff x_i = 0 \right. \right\}\subset \RR^n.\]
(equivalently, $z_i \ne 0 \iff x_i = 1$).
If $x \ne y \in \B_n$ then $Z_x \cap Z_y = \phi$.  If $x < y$ then $Z_x \subset \overline{Z_y}$.
The strata form a decomposition,
\[ \bigcup_{x \in \B_n} Z_x = \RR^n.\]

If $f:\B'_n \to \{0,1\}$ is a ``punctured'' Boolean function\footnote{ cf.~\S \ref{subsec-definitions}.
  Equivalently we could require $f(\ul{\bf 0}) = 0$} define the following ``arrangement"
of linear spaces and complements,
\[ Z(f) \defeq  \bigcup_{x\in \B'_n} \left\{ Z_x|\ f(x) = 1 \right\}.\]

As described in the introduction,  define\footnote{The complexification of this stratification 
$\CC^n = \bigcup_{x \in \B_n} Z_x^{\CC}$ is defined by
replacing  $\RR^n$ with $\CC^n$ in the definition.
It seems likely that $\gamma(f) = \gamma^{\CC}(f)= b^*(Z^{\CC}(f))$ for any Boolean function $f$. A similar thing happens for complements of linear spaces, see \cite{SMT} \S III Cor.~1.4.}
\begin{equation}\label{eqn-gamma} \gamma(f) 
\defeq b^*(Z(f))= \sum_{i \ge 0} \rank H_i(Z(f);\QQ)\end{equation}

\subsection{A formula for $\gamma$.}  There is a remarkable red corner-like formula for the subspaces model.  
Recall from \S \ref{subsec-defs2} the standard $n-1$ dimensional simplex,
$\Delta = \Delta^{n-1}$ with one facet $\partial_i\Delta$ for each $i \in [n]$.
 If $J \subset \{ 1, \cdots, n\}$
let $\partial_J\Delta = \bigcup_{j\in J} \partial_j\Delta$  be the union of the facets
corresponding to the subset $J,$ with inclusion mapping 
\[i_J:|\partial_J\Delta|  \to |\Delta|.\]   
As in  \S\ref{subsec-defs2}
we may identify the barycentric subdivision $\Delta'$ of the $n-1$ dimensional simplex with the order complex of
the punctured Boolean lattice $\B'_n$.  If $f:\B'_n \to\{0,1\}$ is a Boolean function 
we may consider the order complex $\Sigma(f)$ to be a subcomplex of
$\Delta'$. To be explicit, each non-zero Boolean vector
$h = (h_1,h_2,\cdots,h_n)$ (with $h_i \in \left\{0,1\right\}$) corresponds to a vertex
$v(h) \in \Delta'$, the barycenter of the face spanned by $h_1e_1, \cdots, h_ne_n$.  Then 
$\Sigma(f)$ is the full subcomplex of $\Delta'$ spanned by vertices $v(h)$ such that $f(h) = 1$. 

\begin{thm}\label{thm-redcornersubspaces}
Let $f:\B'_n \to\{0,1\}$ be a (punctured) Boolean function.  Then
\begin{equation}
\label{eqn-beta2}\gamma(f) = \sum_{J\subset \left\{1,\cdots,n\right\}}
b^*(|\Sigma(f)|, |\Sigma(f)| \cap |\partial_J\Delta|).\end{equation}
\end{thm}

The proof, which appears in \S \ref{subsec-proof-subspaces}, summarizes the developments
of the next few sections.

\quash{
Recall from \S \ref{subsec-cube} the positive facets $\partial_i^+\c$.  For each subset $J \subset [n]$ let
$\partial^+_J\c = \bigcup_{j\in J} \partial^+_j\c$ be the corresponding union of facets.  
(If $J = \phi$ then $\partial^+_J\c = \phi$.  If $J = [n]$ then $\partial^+_J\c = |K(\B'_n)|$.)

\begin{thm}\label{thm-redcornersubspaces}
Let $f:\B'_n \to\{0,1\}$ be a (punctured) Boolean function.  Then
\begin{equation}
\label{eqn-beta2}\gamma(f) = \sum_{J\subset \left\{1,\cdots,n\right\}}\sum_{i\ge 0} 
\rank{ H_i(\Sigma(f), \Sigma(f) \cap \partial^+_J\c)}.\end{equation}
\end{thm}

The proof, which appears in \S \ref{subsec-proof-subspaces}, summarizes the developments
of the next few sections.

} 

\subsection{The simplicial complex $X(f)$}  \label{subsec-X(f)}
The first step in analyzing the subspaces model is the construction
of a simplicial complex $X(f)$ whose geometric realization $|X(f)|$ is homotopy equivalent to $Z(f)$. 
 Fix $n \ge 1$.  
Let 
 \begin{equation}\label{eqn-Dn}
 \DD^n = \left\{ (x_1,\cdots,x_n) \in \RR^n:\ |x_i| \le 1, \text{ for }
1 \le i \le n \right\}\end{equation} denote the $n$ dimensional cube (or cubical disk) and
let $\SS^{n-1} = \partial \DD^n$ denote its boundary, the {\em cubical sphere}.  
Define  the simplicial complex $\underline{\bf S}^{n-1}$ to be the barycentric subdivision of the 
cell complex $\partial \DD^n$.  Equivalently $\underline{\bf S}^{n-1} = \Sigma (\mathbf{F}^0(\DD^n))$
is the order complex of the poset $\mathbf F^0(\DD^n)$ of {\em proper} faces of the cube $\DD^n$.  Under the natural
identification  $|\underline{\bf S}^{n-1}|\to \SS^{n-1}$, the subset $\SS^{n-1}\cap Z(f)$ is a union of interiors of simplices.

\bigskip
\begin{figure}[h!]
\newcommand*{\V}{.2}
\renewcommand*{\H}{.4}
\begin{tikzpicture}[scale=5]

\draw[thick](0,0) -- (1,0);
\draw[thick](0,0) -- (1,0);
\draw[thick](1,0) -- (1,1);
\draw[thick](0,1) -- (1,1);
\draw[thick](0,0) -- (0,1);
\draw[thick](0,1) -- (\H,1+\V);
\draw[thick] (1,1) -- (1+\H,1 +\V);
\draw[thick]((1,0) -- (1+\H,\V);
\draw[thick] (1+\H,\V) -- (1+\H,1+\V);
\draw[thick] (\H,1+\V) -- (1+\H,1+\V);

\draw[red] (.5,0) -- (.5,1);\draw[red](0,.5) -- (1,.5);
\draw[red] (.5,1) --(.5+\H, 1+\V);
\draw[red](1,.5) -- (1+\H, .5+\V);
\draw[red] (1+\H/2, \V/2) --(1+\H/2, 1+\V/2);
\draw[red] (\H/2,1+\V/2)--(1+\H/2,1+\V/2);
\draw[dotted](0,0) -- (1,1);
\draw[dotted] (1,0) -- (0,1);

\draw[dotted] (1,0) -- (1+\H,1+\V);
\draw[dotted] (1+\H,\V) -- (1,1);
\draw[dotted](0,1) -- (1+\H,1+\V);
\draw[dotted] (1,1) -- (\H,1+\V);
\end{tikzpicture}\caption{The cubical sphere $\SS^2$}
\vspace{9mm}
\end{figure}

The vertices of $ \underline{\bf S}^{n-1}$ are the barycenters of proper faces of $\DD^n$, that is, the points
$x=(x_1,\cdots,x_n)$ such that each $x_i \in \{ -1, 0, 1\}$.  The partial ordering defined by 
 \[x \le y \iff  |x_i| \le |y_i|\quad  \forall i\in [n]\]
 is the {\em reverse} of the partial ordering of the proper faces of $\DD^n$.  For any chain of vertices
\begin{equation}\label{eqn-ordering}
x=x^{(0)} < x^{(1)}< \cdots < x^{(r)}=y\end{equation}  in $\uS^{n-1}$,
their convex hull forms an $r$-dimensional simplex $\sigma$,
with a  standard orientation given by the ordering of its vertices in (\ref{eqn-ordering}).  
The interior $\sigma^o$ is the intersection
\[ \sigma^o = \sigma \cap Z_{|y|}\]
with the stratum $A_{|y|}$  where $|y| = \left(
|y_1|, |y_2|, \cdots, |y_n|\right)$.  The  cubical sphere
$\SS^{n-1}$ is equal to the disjoint union of the interiors of such
simplices $\sigma\in \underline{\bf S}^{n-1}$.   

Now let $f:\B_n' \to \{0,1\}$ be a Boolean function and define
\[ X(f):=  \left\{\sigma \in \underline{\bf S}^{n-1}\left|\
f(|v|)=1 \text{ for every vertex } v \in \sigma \right. \right\} \]
to be the (closed) subcomplex of $\underline{\bf S}^{n-1}$ that is spanned by the
simplices  $\sigma \in\underline{\bf S}^{n-1}$ whose vertices correspond to 
Boolean inputs to which the function $f$ assigns the value $1$.

\begin{prop}\label{prop-homotopy}
For any Boolean function $f:\B'_n \to \{0,1\}$ the following natural inclusions 
are homotopy equivalences:
\[\begin{CD}
 Z(f)@<<<  Z(f) \cap \SS^{n-1} @<{\alpha}<< |X(f)|\end{CD}.\]
\end{prop} 
\begin{proof}  The first inclusion is clearly a homotopy equivalence.
For each $x\in \B'_n$ the connected components of $Z_x \cap \SS^{n-1}$ are cells, that is,
each connected component is homeomorphic to the unit disk in some Euclidean space.  Moreover,
these cells constitute a {\em regular} cell decomposition of $\SS^{n-1}$, meaning that the
closure of each cell is a union of (lower dimensional) cells.  The simplices of
$\underline{\bf S}^{n-1}$ constitute the first barycentric subdivision of this cell complex.
The inclusion $\alpha$ is a homotopy equivalence as a consequence of the following more general statement.
Let $\underline{\bf S}$ be a regular cell complex and let $R \subset | \underline{\bf S}|$ be
a union of interiors of cells.  Let $X \subset R$ be the  union of the set of (closed) simplices 
$\sigma \in \underline{\bf S}'$ in the first barycentric subdivision of $\underline{\bf S}$ such
that $\sigma \subset R$.  Then $X$ is a deformation retraction of $R$, hence
the inclusion $X \to R$ is a homotopy equivalence.
\end{proof}



\subsection{Group action}\label{subsec-groupG}
The group $G= \{\pm 1\}^n$  acts on $\RR^n$ by co\"ordinate-wise
multiplication.  If $f:\B'_n \to \{0,1\}$ is a Boolean function then the group $G$  preserves
$X(f)\subset \SS^{n-1}$.
Hence, it induces a decomposition
\begin{equation}\label{eqn-homology-decomposition}
H_i(|X(f)|;\mathbb Q) = \bigoplus_{\chi \in \widehat{G}} H_i(|X(f)|)_{\chi}\end{equation}
into isotypical pieces, indexed by the characters $\widehat{G}$ of $G,$ 
which in turn, correspond to subsets $J \subset [n]$ with 
$ \chi_J(s_1,\cdots,s_n) = \prod_{j\in J}s_j \in \left\{ \pm1 \right\}$.  In 
\S \ref{subsec-proof-subspaces} we prove the following:
\begin{prop}\label{prop-sums-correspond}
The sum in equation (\ref{eqn-beta2}) corresponding to a subset $J \subset [n]$
is precisely the  isotypical component of $\gamma(f)$ corresponding to the character $\chi_J$.
\end{prop}

\section{Reduction to the order complex $\Sigma(f)$}\label{sec-reduction}
\subsection{Two ways}\label{subsec-pi}
The quotient of $\RR^n$ by the action of the group $G$ (of \S \ref{subsec-groupG}) may be described 
as the map  $\theta:\RR^n \to (\RR_{\ge 0})^n = \RR^n/G$ by
$\theta(x_1, x_2, \cdots, x_n)) = (|x_1|, |x_2|, \cdots, |x_n|)$.
 As in \S \ref{sec-subspaces}  let $\SS^{n-1}$ be the cubical $n-1$ dimensional sphere with its triangulation
${\uS}^{n-1} $. Then 
\[\theta(\SS^{n-1}) = \SS^{n-1}/G = \partial^+\c\cong \left|\Delta'\right|\]  
using the simplicial isomorphism
of equation (\ref{eqn-simplicial-isomorphisms}).  The composition 
\[\pi:\mathbf{S}^{n-1} \to \Delta' \]
is the simplicial map\footnote{The mapping $\pi$ of equation (\ref{eqn-pi-fraction}) is {\em pseudo-linear}:
although its restriction to each simplex fails to be linear, it is homotopic, by a simplex preserving homotopy,
to the piecewise linear map that agrees with $\pi$ on each vertex.  For the purposes of this paper there
is no harm in pretending these are the same map.} 
whose value on each vertex $x = (x_1,x_2,\cdots,x_n)$ is given by
\begin{equation}\label{eqn-pi-fraction} \pi(x_1,x_2,\cdots,x_n) = \frac{1}{\Sigma |x_i|}(|x_1|,|x_2|,\cdots,|x_n|).\end{equation}
It is induced from an {\em order reversing} map of posets $\mathbf F^0(\DD^n) \to \mathbf F(\Delta),$ or equivalently, an order preserving map of posets 
\[\bar\pi:\mathbf F^0(\DD^n)^{op} =\mathbf F^0(\Diamond^n) \to \mathbf F(\Delta)=\B'_n\] from
the poset of proper faces $\mathbf F^0(\Diamond^n)$ of the hyperoctahedron to the poset of nonempty faces of $\Delta$.
If $f:\B'_n \to \{0,1\}$ is a Boolean function then $X(f)$ 
is the order complex of the poset $\bar\pi^{-1}(\B(f))$ and
\[ \pi: X(f) \to \Sigma(f).\]
\smallskip

There are two ways to compute $H^*(|X(f)|)$ in terms of data on $|\Sigma(f)|$.  One method is to use sheaf theory, see \S \ref{sec-sheaves} below.
The other is based on poset topology and the following fiber formula
 of \cite{BWW}.

\begin{prop}\label{posetfibers}
Let $g:P \to Q$ be a poset map such that $Q$ is connected and for all $q \in Q$ the order complex 
$\Sigma(g^{-1}(Q_{\le q}))$ is $\ell(g^{-1}(Q_{<q}))$-connected. Then there is a homotopy equivalence,
$$ \Sigma(P) \simeq \Sigma(Q) \vee \,\bigvee_{q \in Q}
\left(\Sigma(g^{-1}(Q_{\le q})) * \Sigma(Q_{> q})\right).
$$\end{prop}

The wedge on the right side represents a quotient of the disjoint union, where an arbitrarily chosen point
of $\Sigma(g^{-1}(Q_{\le q})$ is chosen (one for each $q\in Q$) and identified with some (arbitrarily chosen)
point of $\Sigma(Q)$.

There is a corresponding result for homology under a corresponding acyclicity assumption which, in the case of
the map $\pi:X(f) \to \Sigma(f)$ says that for each vertex $v$ in $\Sigma(f)$,
 the (reduced) homology of $\pi^{-1}(\Sigma(f)_{\le v})$ vanishes except
possibly in the top dimension.  With this assumption the following formula holds for the Betti sum of reduced homology:

\begin{equation}\label{f1}
\gamma(f) = b^*(X(f)) = b^*(\Sigma(f)) + \sum_{v \in \B(f)} b^*(\pi^{-1}(\Sigma(f)_{\le v})). b^*(\Sigma(f)_{>v}).
\end{equation}
Applications will be given in section \ref{sec-computations}.

\section{Sheaves on the simplex $\Delta$}\label{sec-sheaves}
\subsection{ Sheaves} As in \S \ref{sec-reduction} and equation (\ref{eqn-pi-fraction}) we consider the simplicial map
$\pi: \ul{\bf S}^{n-1} \to \Delta' \cong \ul{\bf S}^{n-1}/G$.
 Denote the constant sheaf on $\SS^{n-1}$ by  $\underline{\mathbb Q}_S$.
 Then its direct image $\pi_*(\underline{\mathbb Q}_{\SS})$ is a
constructible sheaf\footnote{meaning that it is constant on the interior of each simplex of $|\Delta'|.$  
Such a sheaf is then determined by the data of a {\em simplicial sheaf}, that is, a contravariant functor from the category
$\Delta'$ (its faces and inclusions) to the category of Abelian groups.   }
 on $|\Delta'|,$ on which the group $G$ acts.  So it also splits
up as a direct sum
\begin{equation}\label{eqn-sheaf-sum}
 \pi_*(\underline{\mathbb Q}_{\SS}) \cong \bigoplus_{\chi \in \widehat{G}} \mathbf A_{\chi}
\end{equation}
of constructible (or simplicial) sheaves on $\Delta'.$ (This is a general fact about the action of a finite
group on  objects in a $\QQ$-linear  semisimple category.)

\begin{lem}\label{lem-sheaftospace}
  Let $Y \subset |\Delta|$ be a closed subset and let $X = \pi^{-1}(Y).$  
Then the decomposition (\ref{eqn-homology-decomposition}) of the homology of $X$ is 
induced by the ``universal'' decomposition (\ref{eqn-sheaf-sum}) of the sheaf
$\pi_*(\underline{\mathbb Q}_{\SS})$, that is, for any $\chi \in \widehat{G},$ 
\[ H^i(X)_{\chi} \cong H^i(Y; \mathbf A_{\chi}).\]
\end{lem}

\proof We have a Cartesian square of inclusions and projections,
\begin{diagram}[size=2em]
 X & \rTo^{\omega}& \SS^{n-1} \\
\dTo_{\pi} & &\dTo_{\pi} \\
Y & \rTo_{\iota} &|\Delta|
\end{diagram}
Consequently there is a $G$-equivariant isomorphism\footnote{Strictly speaking these 
morphisms and isomorphisms should be interpreted in
the bounded constructible derived category of sheaves on $\Delta|$ however
in this case, the derived functor $R\pi_*$ coincides with $\pi_*.$}
of sheaves
\[\bigoplus_{\chi} \mathbf A_{\chi}|Y \cong \iota^*\pi_*(\underline{\mathbb Q}_{\SS})
\cong \pi_*\omega^*(\underline{\mathbb Q}_{\SS}) \cong
\pi_*\omega^*\pi^*(\underline{\mathbb Q}_{\Delta})
\cong \pi_* \pi^* \iota^* (\underline{\mathbb Q}_{\Delta})
\cong \pi_*(\underline{\mathbb Q}_X)\]
by equation (\ref{eqn-sheaf-sum}).
Using (\ref{eqn-homology-decomposition}), this gives a $G$-equivariant isomorphism
\[ \bigoplus_{\chi\in \widehat{G}}
H^i(Y; \mathbf A_{\chi}) \cong \bigoplus_{\chi \in \widehat{G}}
H^i(X; \mathbb Q)_{\chi}.\]
The result now follows from Schur's lemma.  \qed

It is possible to explicitly describe the sheaves $\mathbf A_{\chi}$ in our case.  
Let $J\subset \left\{1,2,\cdots,n\right\}$ be a subset.  Let $\partial_J\Delta = 
\bigcup_{j\in J} \partial_j\Delta$  be the union of the codimension one faces 
corresponding to the subset $J,$ with inclusion mapping $i_J:|\partial_J\Delta| 
\to |\Delta|.$  Let $(i_J)_* \underline{\mathbb Q}_{\partial_J\Delta}$
denote the constant sheaf on $|\partial_J\Delta|,$ viewed as a sheaf on $|\Delta|.$

\begin{prop}\label{prop-sheaf-sum}
In the decomposition (\ref{eqn-sheaf-sum}), the sheaf $\mathbf A_J$
corresponding to the character $\chi_J$ is
\[ \mathbf A_J = \ker\left( \underline{\mathbb Q}_{\Delta} \to 
(i_J)_*\underline{\mathbb Q}_{\partial_J\Delta}  \right).\]
\end{prop}
The proof will be given in Appendix \ref{appendix1}.

In other words, the sheaf $\mathbf A_J$ is the constant sheaf on
$|\Delta| - |\partial_J\Delta| $ and it is zero on the (closed) face $|\partial_J\Delta|.$
The sheaf $\mathbf A_J$ may also be described as the sheaf on $|\Delta|$ that is obtained
from the constant sheaf on $|\Delta| - |\partial_J\Delta|$ by extension by zero.
In more down-to-earth terms, if $Y \subset |\Delta|$ is a closed subspace then
its cohomology is given by the relative cohomology group
\begin{equation}\label{eqn-relative-cohomology}
H^i(Y; \mathbf A_{J}) = H^i(Y, Y \cap |\partial_J\Delta|; \mathbb Q). \end{equation}

\subsection{Proof of Theorem \ref{thm-redcornersubspaces} and Proposition \ref{prop-sums-correspond}}
\label{subsec-proof-subspaces}
By Proposition (\ref{prop-homotopy}), the quantity $\gamma(f)$ is the sum of Betti numbers of 
the subset $|X(f)| \subset \SS^{n-1}$ of the ``cubical sphere".
The projection $\pi:\SS^{n-1} \to |\Delta'|$  takes $X(f)$ to the order complex
$\Sigma(f)$ but its cohomology breaks up under the action of the group $G\cong (\ZZ/(2)^n$, 
which, according to Lemma \ref{lem-sheaftospace}  is a reflection of the fact 
that the push forward of the constant sheaf breaks up under this action into sheaves $A_{\chi}$.
These sheaves $A_{\chi}$ are explicitly described in Proposition (\ref{prop-sheaf-sum}), and the
effect on cohomology is described in equation (\ref{eqn-relative-cohomology}) 
(with $Y = |\Sigma(f)|)$ as  giving the relative cohomology groups
that appear in Theorem \ref{thm-redcornersubspaces}. That is: 
\[ H^i(Z(f)) \cong H^i(|X(f)|) \cong \bigoplus_{\chi}H^i(|\Sigma(f)|; A_{\chi}) \cong \bigoplus _{J\subset [n]}
H^i(|\Sigma(f)|, |\Sigma(f)| \cap |\partial_I\Delta|).\]
It also follows that $\beta(f) \le \gamma(f)$ since $f({\bf\ul 0})=0$ and for each red corner,
 $R_v \cap |\Sigma(f)|= |\partial_J\Delta| \cap |\Sigma(f)|$, where $J = [n]-\{v\}$.\qed

 \section{Proof of Theorem \ref{thm-intro}}  
 \label{sec-proof-subspaces}

A proof 
for the case of the Red Corner model was given 
in Section \ref{sec-RedCornerProof}.

 \subsection{Proof of Theorem \ref{thm-intro} for the subspaces model} 

Part ( 1) of Theorem \ref{thm-intro} is clear from the definition of $Z(f)$.  
Part (2) will be proven below.

Part (3), the K\" unneth theorem
applied to $Z(f \wedge g) \cong Z(f) \times Z(g)$ implies $\gamma(f \wedge g) = \gamma(f)\gamma(g)$ and
the same result for $\gamma(f \vee g)$ follows from this and part (2), provided neither
$f$ nor $g$ is identically $0$ or $1$.  Part (4) is clear since $Z(f) \cong Z(g) \times \RR$.
Part (5) follows from parts (2) and (3) and induction by considering the top gate since each branch involves
different input variables.  The induction begins with $n=2$, noting that $\gamma(\text{OR}) = 2$ and 
$\gamma(\text{AND}) = 4$.

The proof of part (2) requires some technical results from the Appendix \ref{sec-Verdier}.   
Given Boolean functions $f:\B'_n \to \{0,1\}$ and $g = \sim f$ 
 then $\Sigma(f)$ and $\Sigma(g)$ are (geometric realizations of) {\em supplementary} subcomplexes of the simplex 
 $\Delta'$ (cf. Lemma \ref{lem-P1}).  As above, decompose 
 \begin{equation}\label{eqn-Jsum}
 H^*(|X(f)|) = \bigoplus_{I\subset[n]}
 H^*(|\Sigma(f)|, A_I)\qquad H^*(|X(g)|) = \bigoplus_{J \subset [n]} H^*(|\Sigma(g)|, A_J).\end{equation} 
 We will compare these sums, term by term.
 If $I,J \subset [n]$ are complementary subsets then  by Lemma \ref{lem-Verdier-dual}
 the sheaves $A_{I}=A_{\chi_I}$ and  $A_J=A_{\chi_J}$ are (Verdier) dual to each other (up to a shift).

First consider the case when $I,J \ne \phi$.   Then $H^r(|\Delta|, \mathbf{A}_J) = 0$ for all $r$, 
because the relative cohomology of the simplex $|\Delta^n|$ modulo any nonempty connected contractible
 subset vanishes in all degrees.  So if $I,J$ are proper subsets of $[n]$ and complementary, then
 the sheaves $\mathbf {A}_I, \mathbf{A}_J$ are acyclic and Verdier dual (up to a shift).
 Taking $\mathbf A\b = \mathbf{A}_J$ and 
$\mathbf T\b = \mathbf{A}_I[n-1]$ 
in Proposition \ref{prop-Alexander}  gives cohomology isomorphisms
 \[ b^*(|\Sigma(f)|, \mathbf{A}_I) = b^*(|\Sigma(g)|, \mathbf{A}_J).\]

 If $I = \phi$ then $A_I = \QQ_{\Delta}$ and $J = [n]$ so  $A_J $ is the constant sheaf on the 
 interior of $|\Delta|$, extended by zero  (which is denoted $\QQ_!$ in \S \ref{subsec-AJ}).  So we must prove that
 \begin{equation} \label{eqn-fg-equality} b^*(|\Sigma(f)|) + b^*(|\Sigma(f)|, |\Sigma(f)| \cap |\partial\Delta|) =
 b^*(|\Sigma(g)|) + b^*(|\Sigma(g)|, |\Sigma(g)| \cap |\partial\Delta|).\end{equation}
 We may suppose that $f(\ul{\mathbf 1}) = 1$ (otherwise $g(\ul{\mathbf 1}) = 1$).  Let $\widehat{f}$ be the Boolean function
 $\widehat{f}(\ul{\mathbf 1}) = 0$ and $\widehat{f}(v) = f(v)$ for $v \ne \ul{\mathbf 1}$.   Then 
 $|\Sigma(f)| \cap |\partial\Delta| = |\Sigma (\widehat{f})|$
 and $|\Sigma(f)|$ is the cone over $|\Sigma(\widehat{f})|$.  It follows that $b^*(|\Sigma(f)|) = 1$ and
 $b^*(|\Sigma(f)|, |\Sigma(f)| \cap |\partial \Delta|) = b^*(|\Sigma(\widehat{f})|) -1$.  (In fact, it is the total Betti number
 of the reduced cohomology of $|\Sigma(\widehat{f})|$.)
 
 On the other hand, $b^*(|\Sigma(g)|, |\Sigma(g)|\cap |\partial \Delta|) = 0$ since $\Sigma(g) \subset \partial\Delta$.
 If $g$ is {\em not identically zero} then the complexes $|\Sigma(\widehat{f})|$ and $|\Sigma(g)|$ are Alexander dual 
 within $\SS^{n-1}$:  the complement of one deformation retracts to the other.  Thus,
 \[ b^*(|\Sigma(g)|) = b^*(|\partial\Delta| - |\Sigma(\widehat{f})|) = b^*(|\Sigma(\widehat{f})|))\]
 by Alexander duality, which states that the reduced homology of $|\Sigma(\widehat{f}))|$ equals the reduced
 cohomology of $\SS^{n-1} - |\Sigma(\widehat{f})|$ so they have the same total Betti sum.  Thus we have shown
 that both sides of equation (\ref{eqn-fg-equality}) are equal to $b^*(|\Sigma(\widehat{f})|)$. 

The theorem omits the cases when $f\equiv 0 \text{ or } 1$, since $\gamma(f) = 2$ if $f \equiv 1$ but
$\gamma(\sim f) = 0$. \qed
 

\section{Examples}\label{sec-computations}

In this section we discuss estimates for the circuit complexity of
some of the most basic Boolean functions, PARITY (see \S \ref{subsec-parity}), 
MAJORITY and THRESHOLD (see \S \ref{subsec-threshold}), and the role
of conjecture (\ref{eqn-conjecture}).

\begin{prop}  
If the conjecture (\ref{eqn-conjecture})  is true, then
\begin{align*}
CC(f^{\mathtt{PARITY}}_{n}) &\ge \frac{1}{5}n^2,\\
CC(f^{\mathtt{MAJORITY}}_{n}) &\ge 2\sqrt{2}\ n, \\
CC(f^{\mathtt{THRESHOLD}}_{n,[2n/3]}) &\ge 3\ n. 
\end{align*}
\end{prop}
Knowing that random Boolean functions have exponential circuit
complexity, these bounds may seem disappointingly modest. However, 
these bounds have  the right order of magnitude, comparable to the best 
bounds known in the literature, see e.g. \cite{We}.

\subsection{The parity function}\label{subsec-parity}

The function  $f^{\mathtt{PARITY}}_{n}(s_1,s_2,\cdots,s_n)$ 
equals $1$ if $s_i=1$
for an even number of the $n$ inputs, and equals 0 otherwise.
The quadratic lower bound given above for 
$\mathtt{PARITY}$  can with a little more care be improved to  $\frac{1}{4}n^2$.
In fact (\cite{Krapchenko1,Krapchenko2}), the circuit complexity of 
parity is known to be exactly $n^2$, so our estimate is not bad.

To prove the bound for the parity function, we need to 
use some properties of {\em Euler numbers} $E_n$, defined by
$\sum \frac{E_n}{n!} x^n = \sec(x) + tan(x)$. They
have a well-known combinatorial  interpretation (see \cite{Alternating}), namely,
$$E_n =  \text{the number of alternating permutations in the symmetric group}
\ {\mathbb S_n}.$$
A permutation $\pi: [n] \rightarrow [n]$  is said to be {\em alternating}, if
$$\pi(1) < \pi(2) > \pi(3)  < \pi)4) > \pi(5)< \cdots $$

For the Boolean lattice $\B_n$ and any subset
$J\subseteq \{1,2, \ldots ,n\}$, consider the {\em rank-selected subposet} 
$\B_J\defeq \{ v\in \B_n  \mid \rank(v) \in J \}$. Its order complex
has the following enumerative 
properties, quoted from Theorem 10 in \cite{CohenMacaulay} and 
Corollary 3.13.2 in \cite{EC1}. The {\em descent set} of $\pi \in {\mathbb S_n} $ is 
$D(\pi) \defeq\{j \mid  \pi(j) > \pi (j+1)  \}$.

\begin{prop}{} \label{rankselect}

\begin{enumerate}
\item The simplicial complex $\Sigma(\B_J)$ is Cohen-Macaulay\footnote{A complex is said to be 
{\em Cohen-Macaulay} if its homology vanishes below 
the top dimension, both for the complex itself and for the links of all its faces.
For such a complex, the sum of Betti numbers equals the top Betti number,
which in turn equals the reduced Euler characteristic, up to sign.}. 
\item The sum of Betti numbers of $\Sigma(\B_J)$
 is equal to the number of permutations of
${\mathbb S_n}$  with descent set $J$. \\
\end{enumerate}
\end{prop}
Specialized to the case of interest  $J=D(\pi)=\{2,4,6, \ldots \}$, 
this implies that the sum af Betti numbers of
$\Sigma(\B_J) = \Sigma(f^{\mathtt{PARITY}}_{n} )$ equals $E_n$.
Thus we have 
\begin{equation} 
\gamma (f^{\mathtt{PARITY}}_{n} ) \ge
 b^*(\Sigma(f^{\mathtt{PARITY}}_{n} )) = E_n:
\end{equation}
where the inequality comes from Proposition (\ref{f1}).


From \cite{Analytic} or \cite{Alternating} we have the following estimate for Euler numbers:
\[ E_{2n} > 8 \sqrt{\frac{n}{\pi}} \left( \frac{4n}{\pi e}\right)^{2n} .\]
Thus, for even  numbers $n$,

\[ E_{n} >  \frac{8}{\sqrt{2\pi}}{\sqrt{n}} \left( \frac{2n}{\pi e}\right)^{n}
\ge
   \left( \frac{n}{c}\right)^{n} \ge \left( \frac{n}{5}\right)^{n}, \]
   where $c=\frac{e\pi}{2} \cong 4.27$. 
  The algebraic-geometric mean inequality (\ref{eqn-geometric-mean}) then gives 
\[ CC(f^{\mathtt{PARITY}}_{n})
\ge n\cdot (\gamma{(f^{\mathtt{PARITY}}_{n}))^{1/n} \ge n\cdot \frac{n}{5}},
\]
as was to be shown.
\smallskip

{\bf Remark.}
Formula (\ref{f1}) applies and yields the following recursive 
description of the  $\gamma$ function of parity (details omitted): 

$$\gamma(f^{\mathtt{PARITY}}_{n}) =B_n+1,$$ 
where
$$B_n = E_n +\left\{
\begin{array}{lll}
\sum_{k=1}^{(n-1)/2} \binom{n}{2k} E_{n-2k} B_{2k}&& \mbox{if $n$ is odd}  \\
&\\
\sum_{k=1}^{n/2} \binom{n}{2k-1} E_{n-2k+1} B_{2k-1} && \mbox{if $n$ is even}
\end{array}
 \right.
$$
We compute from this:
\begin{center}
\begin{tabular}{||l||l|l|l|l|l|l||}
\hline & $n=1$ & $n=2$ & $n=3$ & $n=4$ & $n=5$&$n=6$\\
\hline $E_n$& 1& 1 & 2 & 5 & 16 & 61  \\
\hline $B_n$& 1 &  3 & 11 & 57 & 361 & 2763  \\
\hline 
\end{tabular}
\end{center}
The numbers $B_n = \gamma((f^{\mathtt{PARITY}}_{n}) -1 $ 
appear in the online encyclopedia of integer sequences, 
where they are identified as the exponential 
generating function of $1/(1-\tan(x))$.

\subsection{Monotone functions.}
Suppose $f:\B_n \to \{0,1\}$ is monotone, meaning that $v \le w\implies f(v) \le f(w)$.  Equivalently,
the function $f$ is computable by a circuit with no NOT gates.  For the purposes of this section
only, we temporarily consider $\B(f) = f^{-1}(0)$ rather than $f^{-1}(1)$, because it is easier to keep track of the relations this way. 
If $v \in \B(f)$ then $f(v) = 0$ so every $w\in \B_n$ with $w \le v$ is
also in $\B(f)$.  Thus, $\B(f)_{\le v}$ is the poset of faces of a simplex of dimension $|v|-1$.  Therefore 
$\B(f)$ can be interpreted as the face poset of an abstract
simplicial complex $K_f$ with one simplex $\langle v \rangle$ for each $v \in \B(f)$.  Moreover, the
pre-image $\bar\pi^{-1}(\B(f)_{\le v})$ is the full poset $\mathbf{F}^0(\Diamond^{|v|})$ of proper faces
of the $|v|$-dimensional hyperoctahedron (cf.~\S \ref{subsec-pi}), whose order complex is topologically
a sphere of dimension $|v|-1$.  Also, $\B(f)_{>v}$ is the face poset of the link of the simplex $\langle v\rangle$ in the
simplicial complex $K_f$.
Thus, for monotone $f$ formula (\ref{f1}) specializes to
\begin{equation}\label{f2}
\gamma(f) = b^*(K_f) + \sum_{v\in \B(f)} b^*(\mathrm{link}_{K_f}(\langle v \rangle)).
\end{equation}
 
In this case $\Sigma(f)$ can be interpreted as (the barycentric subdivision of) the simplicial complex $K_f$,
and $X_f$ has a simplicial interpretation as the $2$-inflation of $K_f$ (in the language of
\cite{BWW}) so that formula (\ref{f2}) is a consequence of
 \cite[Cor.6.3]{BWW}. 

\subsection{Symmetric functions}
A Boolean function $f(s_1,s_2,\cdots,s_n)$ is said to be {\em symmetric} 
or {\em cardinality-based},
if its value depends only on the number of input variables $s_i$ that are
equal to one. For $J\subseteq \{1,2, \ldots ,n\}$, we
denote by $f_J$ the function that takes value one
precisely if the number of input ones is a number belonging to $J$.
Threshold functions and parity functions are of this kind.
The formula (\ref{f1}) is applicable to all symmetric functions, because
rank-selection preserves the Cohen-Macaulay property (see \cite{CohenMacaulay}). 

  If $u \le v \in \B_n$ let $\Sigma(u,v)$ be the order complex of
the open interval $(u,v)$.  For any simplicial complex $\Sigma$ let $H_{\rm top}(\Sigma) = H_{\dim(\Sigma)}(\Sigma)$.
For the red corner model we have:

\begin{prop}
Suppose $f$ is symmetric and $f(\ul{\mathbf 0}) = f(\ul{\mathbf 1}) = 1$.  Then
\begin{equation}\label{eqn-symmetric} \beta(f) \ge \sum_{ \{v| f(v) = 0 \} }
\mathrm{rank } H_{\rm top}(\Sigma(\ul{\mathbf 0},v) \cap \Sigma(f))\cdot
\mathrm{rank } H_{\rm top}(\Sigma(v, \ul{\mathbf 1}) \cap \Sigma(f)).\end{equation}
\end{prop}
\begin{proof}
By Proposition \ref{prop-Gammav} it suffices to show that the quantity on the right side of
equation (\ref{eqn-symmetric}) is bounded by 
\[ \sum_{i \ge -1} \mathrm{rank} \tilde{H}_i(\Gamma_v(f)).\]
Let $m_1$ be a maximal chain in $\Sigma(\ul{\mathbf 0}, v) \cap \Sigma(f)$ and $m_2$ a maximal chain in
$\Sigma(v, \ul{\mathbf 1})\cap \Sigma(f)$. Then the concatenation $m_1*m_2$ is a maximal chain in $\Gamma_v(f)$. 
The totality of all such concatenated chains $m_1 *m_2$ are the facets of
the simplicial join of the order complexes of $\Sigma(\ul{\mathbf 0}, v) \cap \Sigma (f)$
 and $\Sigma(v, \ul{\mathbf 1}) \cap \Sigma(f)$.
These complexes are homotopy-wedges-of spheres in the top dimension, by a
general fact about rank-selection in Boolean lattices.

Hence, the cycle space of these concatenated chains has rank equal to the
product of the ranks of $H_{\mathrm{top}}(\Sigma(\ul{\mathbf 0},v) \cap \Sigma(f))$ and
$H_{\mathrm{top}}(\Sigma(v, \ul{\mathbf 1}) \cap \Sigma(f))$.  In the top dimension
there are no boundaries so these concatenated cycles contribute independently
to the top homology $H_{\mathrm{top}}(\Gamma_v(f))$. 
\end{proof}

\subsection{Threshold functions}\label{subsec-threshold}

The function 
$ f^{\mathtt{THRESHOLD}}_{n,k}(s_1,s_2,\cdots,s_n)$
 equals 1 if $s_i=1$
for at least $k+1$ of the $n$ inputs, and equals 0 otherwise.
The majority function ($k = n/2$) is a special case.
The class of threshold functions equals the intersection of
the classes of monotone and symmetric function.
The following can be said about the 
problem of maximizing  $\gamma(f_J)$ for threshold functions.

\begin{prop}
\begin{enumerate}
\item The total Betti number of a threshold function is
\[ \gamma(f^{\mathtt{THRESHOLD}}_{n,k}) 
= \sum_{j=0}^k \binom{n}{j}\binom{n-j-1}{k-j}.\]
\item Among all symmetric functions  $f_J$, the total Betti number $\gamma(f_J)$ is maximized by $f^{\mathtt{PARITY}}_{n} $.
\item Among all monotone symmetric functions  $f_J$, the total Betti number
$\gamma(f_J)$ is maximized by 
$f^{\mathtt{THRESHOLD}}_{n,[2n/3]} $ 
\end{enumerate}
\end{prop}

\begin{proof}
Formula (\ref{f2}) applies and specializes to the formula of part $(1)$.  
Parts $(2)$ and $(3)$ follow from work by M. Readdy on the M\"obius function
of rank-selected subposets of Boolean lattices \cite[\S 3]{Readdy}.
\end{proof}

We can estimate the magnitude
of the total Betti number for these Boolean functions using Stirling's approximation $n!\sim (\frac{n}{e})^n\sqrt{2\pi n}$. For instance,

$$ 
\gamma (f^{\mathtt{THRESHOLD}}_{n,[2n/3]})  \ge 
\binom{n}{[n/3]} \binom{[2n/3]-1}{[n/3])}  \ge 3^n
$$

Thus, we reach the following estimate for the most expensive 
threshold function:
$$CC(f^{\mathtt{THRESHOLD}}_{n,[2n/3]}) \ge 3\ n. $$



\subsection{Graphs}
Boolean functions arise from graphs and higher multigraphs in a very
natural way. 
A graph on $n$  labeled vertices  may be interpreted as an element of the (punctured) Boolean lattice
$ \B_m$ where $m = \binom{n}{2}$, that is, a vector of $1$s (representing an edge) and $0$s 
(representing no edge) with coordinates given by unordered pairs of vertices,
 or in symbols,
\[ \B_m = \{0,1\}^{\binom{V}{2}},\]
where $V$ denotes the set of vertices.
Graphs on $n$ vertices are partially ordered by inclusion.  A Boolean function $f:\B_m \to \{0,1\}$ is a yes/no function defined on the set of all graphs with $n$ vertices. It is sometimes fruitful to think a Boolean function $f$ in this setting 
as a ``graph property", the graphs having the property in question
being the ones mapped by $f$ to $1$. We mention below a couple of examples,
for many more see \cite{Jonsson}.

A $k$-{\em clique} in a graph $G$ is a subset of $k$ vertices in $G$  spanning an induced complete
subgraph.  A $k$-{\em anticlique} (= independent set = stable set) is a collection of $k$ vertices in $G$
with no edges between any of them.  Let
$c_{n,k}$ be the Boolean function on the set $\B_m$ of graphs with $n$ vertices which verifies that a
given graph contains no clique of size $k$.  Let $a_{n,k}$ be the Boolean function on the set of graphs
with $n$ vertices which verifies that a given graph contains no anticlique of size $k$.  Then $a_{n,k}$ and
$1-c_{n,k}$ are monotone functions.  

The corresponding complexes $A_{n,k} = \Sigma(a_{n,k})$ and $C_{n,k} = \Sigma(c_{n,k})$ are
{\em Alexander dual} in the simplicial sense, which implies, (see \cite{Bjorner2}) for all $j$ that
\[ \tilde H_j(A_{n,k}) = H_{\binom{n}{2}-j-3}(C_{n,k}).\]
The maximal faces of $C_{n,k}$ are known from Turan's theorem \cite{Turan, Turan2}.
 For each partition $\tau$ of the vertex set $\{1,2,\cdots,n\}$ into
$k-1$ blocks let $G_{\tau}$ be the complete $(k-1)$-partite graph having all edges that connect vertices in
different blocks but no edges within a block.  Each $G_{\tau}$ (together with all of its sub-graphs) forms a
maximal simplex in $C_{n,k}$ and there are $S(n, k-1)$ (Stirling number of the second kind) such choices $\tau$.  
Assume for simplicity that $k-1$ divides $n$.  Then the number of edges in $G_{\tau}$ is
\[ \dim(C_{n,k}) = \binom{k-1}{2}\left( \frac{n}{k-1}\right)^2 -1 = \frac{n^2(k-2)}{2(k-1)}-1.\]
Homology groups of $C_{n,3}$ (that is, triangle-free graphs) were computed by Jonsson \cite{Jonsson} \S 26.7:

\newcommand{\Z}{\mathbb{Z}}
\begin{center}
\begin{tabular}{|c||c|c|c|c|c|c|c|c|c|c|}
\hline $H_i(C_{n,3})$& $i=2$ & $i=3$ & $i=4$ & $i=5$ & $i=6$ & $i=7$ & $i=8$\\
\hline $n=4$& $\Z^3$& 0 & 0 & 0 & 0 & 0 & 0   \\
\hline $n=5$& 0 & $\Z^5$ & $\Z$ & 0 & 0 & 0 & 0   \\
\hline $n=6$& 0& 0 & $\Z^6$ & $\Z^9 \oplus \Z_2$ & 0 & 0 & 0   \\
\hline $n=7$&0& 0 & 0 & $\Z^7$ & $\Z_2$ & $\Z^{55}$ & $\Z$   \\
\hline
\end{tabular}
\end{center}
\medskip

\section{Further speculation}\label{sec-speculation}
We will describe another source of topological spaces that may be associated to any Boolean function.

Suppose a projective complex algebraic variety $X \subset \CC P^N$ is preserved by the action of an algebraic torus $T = (\CC^{\times})^n$.  Suppose the fixed point set $F = X^T$ consists of finitely many points.  The action of the real torus $(S^1)^n \subset T$ determines a moment map $\mu: \CC P^N \to \RR^n$.  We are interested in the cases when the set $F$ of fixed points contains $2^n$ elements and the image $\mu(X) \subset \RR^n$ is an $n$-dimensional cube, whose vertices are the $2^n$ distinct points $V=\mu(F)$.  

If $\mathcal O\subset X$ is an orbit of $T$ then its closure $\overline{\mathcal O}$ contains various fixed points:  $\overline{\mathcal O}^T \subset F$.  For example, the closure of any $1$ dimensional orbit contains exactly two fixed points.
For any subset $S \subset  V$ of the vertices, define the $T$-invariant closed subvariety $X_S \subset X$ to be the union
of all $T$-orbits $\mathcal O$ with the property that $\mu(\overline{\mathcal O}^T) \subset S$. 

Choose an identification of the vertices $V$ with Boolean vectors.  (There may or may not exist a natural choice for this identification.)  Then a Boolean function $f: \{0,1\}^n \to \{0,1\}$ may be considered a mapping $F:V \to \{0,1\}$ and it determines a subset $S=F^{-1}(1) \subset V$. From the resulting subvariety $X_S \subset X$  we may associate a ``complexity measure"
\[ \alpha_X(f) = b^*(X_S) =  \sum_{i \ge 0} \text{rank } H^i(X_S).\]

One example of such a variety $X$ is the (complex) Lagrangian Grassmannian of $n$-dimensional Lagrangian subspaces of $\CC^{2n}$ (with its usual symplectic $2$-form $\omega$).  It is a nonsingular projective algebraic variety of dimension $n(n+1)/2$ and it admits the action of an $n$ dimensional torus.   The standard basis  $\{e_1, e_2, \cdots, e_n, f_1, f_1, \cdots, f_n\}$
of $\CC^{2n}$ is a symplectic basis, that is, 
$\omega(e_i, f_j) = \delta_{ij}$ and $\omega(e_i, e_j) = \omega(f_i, f_j) = 0$.  The standard torus, consisting of diagonal
symplectic matrices $ \left(\begin{smallmatrix} t & 0 \\ 0 & t^{-1} \end{smallmatrix}\right)$ (with $t \in  (\CC^{\times})^n$)
takes Lagrangian subspaces to Lagrangian subspaces, which defines an action on the Lagrangian Grassmannian.

The fixed points of this torus action are the Lagrangian subspaces defined by the 
multivectors $g_1\wedge g_2\wedge \cdots \wedge g_n$
where each $g_i = e_i \text{ or } f_i$. So there is a natural identification of fixed points with Boolean vectors, with
$(0,0,\cdots, 0)$ corresponding to $e_1\wedge e_2\wedge \cdots \wedge e_n$ and with $(1, 1, \cdots, 1)$ corresponding to
$f_1 \wedge f_2 \wedge \cdots \wedge f_n$ and in fact the moment map may be chosen so that the image 
$\mu(X)\subset \RR^n$ is precisely the standard $n$-cube with sides of unit length.

In summary, this gives a {\em Lagrangian Grassmannian model} for complexity of Boolean functions.  We do not know if
there is a useful relation between this number $\alpha_X(f)$ and the circuit complexity $CC(f)$.

\appendix
\section{Decomposition of $\theta_*(\QQ_{\SS^{n-1}})$}\label{appendix1}
\subsection{}
In this appendix we prove Proposition \ref{prop-sheaf-sum}.
As in \S \ref{subsec-X(f)} let $\DD^n$ denote the $n$ dimensional cube (\ref{eqn-Dn}) with
its boundary, the cubical sphere $\SS^{n-1}$.
Dividing by reflections gives a covering $\theta:\DD^n \to \c=[0,1]^n$ with group $G=(\ZZ/(2))^n$.

Pushing forward the constant sheaf gives a sheaf $T_n=\theta_*(\QQ_{\DD^n})$ which
decomposes under the action of $G$ into isotypical components that we now determine.

First consider the case $n=1$ and write $\DD = \DD^1$, $I = I^1=[0,1]$ and $T=T_1$
so that $T$ has stalk $\QQ \oplus \QQ$ at
points in $(0,1]$ and stalk $\QQ$ over $\{0\}$.  The group $G= \ZZ/(2)$ acts 
on $T$ and the subsheaf $T(1)$ corresponding to the trivial character $(1)$ 
is the constant sheaf $\QQ_{[0,1]}$ with inclusion $T(1) \to T$
given by $x \mapsto (x,x) \in \QQ \oplus \QQ$ over points in $(0,1]$. 
 The subsheaf $T(-1)$ corresponding to the nontrivial character $(-1)$ has stalk 
 zero at $\{0\}$ and stalk $\QQ$ at other points with inclusion
$T(-1) \to T$ given by $y \mapsto (y,-y)$ so that
\[ T = \theta_*(\QQ_{\DD}) = T(1) \oplus T(-1).\]
Thus, $T(1) = \QQ_{[0,1]}$ is the constant sheaf and $T(-1) = \ker\left(\QQ_{[0,1]} \to \QQ_{\{0\}}\right)$
or equivalently, it is the extension by zero of the constant sheaf $\QQ$ on $(0,1]$.  See also Lemma \ref{lem-circle}.
The map $[-1,1] \to [0,1]$  is illustrated in Figure \ref{fig-theta}.

\begin{figure}[h!]
\begin{tikzpicture}[scale = 2]
\draw[fill](0,0) circle [radius = .025];
\draw [very thick] (0,0) -- (2,0);
\draw[fill] (2,0) circle [radius = .025];
\draw [very thick] (0,0) -- (2,1);
\draw [fill] (2,1) circle [radius = .025];
\node[right] at (2,1) {\ $-1$};
\node[right] at (2,0) {\ $1$};
\node[left] at (0,0) { $0$\ };

\draw[very thick, -latex] (1,-.2) -- (1,-.7); 
\node[] at (1.2, -.4) {$\theta$};

\draw[fill] (0,-1) circle [radius = .025];
\draw[fill](2,-1) circle [radius = .025];
\draw[very thick] (0,-1) -- (2,-1);

\node[above] at (0,-1) {$\mathbb Q$};
\node[below] at (1,-1) {$\mathbb Q \oplus \mathbb Q$};
\node[above] at (2,-1) {$\mathbb Q \oplus \mathbb Q$};
\end{tikzpicture}
\caption{The mapping $\theta:\DD^1 \to [0,1]$}\label{fig-theta}
\end{figure}

\subsection{} In general,
for each $1 \le j \le n$ let $p_j:\c \to [0,1]$ be projection $(x_1,\cdots,x_n) \mapsto x_j$ to the j-th coordinate and let
$\partial_j^- \c$ be the facet where $x_j = 0$.  
The map $\theta:\DD^n \to 
\c$ is a product of copies of the $n=1$ case and $T_n = \bigotimes_{j=1}^n p_j^*(T)$.  
If $J \subset [n]$ let $\partial_J^- \c = \bigcup_{j\in J} \partial_j^- \c$ with inclusion $i_J:\partial_J^-\c \to \c$.
Then
\[\bigotimes_{j\in J} p_j^*(T(-1)) = \ker\left(\QQ_{\c} \to \QQ_{\partial_J^- \c}\right).\] 
It is the extension by zero (to all of $\c$) of the constant sheaf on $\c-\partial_J^-\c$. 
It follows that  
\[T_n=\theta_*(\QQ_{\DD}) = \bigoplus_{J \subset [n]} T_n(\chi_J) 
= \bigoplus_{J \subset [n]} \ker\left(\QQ_{\c} \to (i_J)_*\QQ_{\partial_j^-\c}\right)
\]

\quash{
\subsection{} 

  As in equation (\ref{eqn-simplicial-isomorphisms}) 
 the identification $|\Sigma(\B_n)| \cong \c$ gives a simplicial isomorphism 
 $\Delta'=|\Sigma(\B'_n)| \cong \partial^+\c$.
 Let $f:\B'_n \to \{0,1\}$ be a (punctured) Boolean function.
Recall from \S \ref{subsec-cube} that the unit cube $\c=[0,1]^n$ has
$2^n$ facets, $\partial_j^{\pm}\c$ where $1 \le i \le n$.  Those containing $\ul{\mathbf 1}$ 
form the {\em positive sphere} $\partial^+\c$
and the reflection mapping $\theta:\DD^n \to \c$ gives a diagram

\begin{equation}\label{eqn-diagram}
\begin{diagram}[height=1.5em]
\DD^n & \rTo^{\theta} & \c & \cong & \ \left| \Sigma(\B_n)\right| \\
\uInto && \uInto^{ j\  } && \uInto \\
\SS^{n-1} &\rTo & \partial^+\c & \cong & \ \left| \Sigma(\B'_n)\right| \\
\uInto && \uInto && \\
X(f) & \rTo & \Sigma(f) &&
\end{diagram}
\end{equation}
}

\subsection{}  \label{subsec-partial_J}
Let us restrict these sheaves to $\partial^+\c $.
If $j:\partial^+\c \to \c$ is the inclusion then 
\[ \theta_*(\QQ_{\SS^{n-1}}) = j^*T_n = \bigoplus_{J\subset[n]} \mathbf A_J \cong
\bigoplus_{J \subset [n]} \ker \left(\QQ_{\partial^+\c} \to (i_J)_*\QQ_{\partial_J^-\c}
\right) \]
Using the simplicial isomorphism $\tau:\partial^+\c \to \Delta'$ we may interpret these as sheaves on $\Delta'$.
 By  (\ref{eqn-simplicial-isomorphisms}) and  (\ref{eqn-facet-intersection}) this isomorphism  takes
 $\partial^+\c \cap \partial^-_J\c$  to $\partial_J\Delta$. Therefore the term $\mathbf A_J$ corresponding to the subset
 $J \subset [n]$ is:
\[ \mathbf A_J= \ker \left( \QQ_{\Delta} \to (i_J)_*\QQ_{\partial_J\Delta}
\right) .\]
which completes the proof of Proposition \ref{prop-sheaf-sum}.  \qed

\section{Alexander-Verdier duality}\label{sec-Verdier}
\subsection{}\label{subsec-supplementary}

Let $X$ be a finite simplicial complex.
Recall that a complex of sheaves $\mathbf A\b$ (of $\QQ$-vector spaces) on $|X|$ is  {\em
cohomologically constructible} (with respect to the simplicial decomposition) if each of its cohomology sheaves
$\mathbf H^i(\mathbf A\b)$ has finite rank and is constant on the interior of each simplex.
If $\mathbf T\b$ is the Verdier dual\footnote{meaning that there is a homomorphism
$\mathbf A\b \otimes \mathbf T\b \to \omega_X$ that induces a quasi-isomorphism
$\mathbf T\b \cong {\rm RHom}(\mathbf A\b, \omega_X)$ where $\omega_X$ is
the dualizing sheaf on $X$.} of $\mathbf A\b$ then $\mathbf T\b$ is also constructible and Verdier
duality gives an isomorphism $H^{-r}(|X|, \mathbf T\b) \cong \Hom(H^r(|X|, \mathbf A\b), \QQ)$
for each integer $r$.

\begin{prop}\label{prop-Alexander}
Let $\mathbf A^{\bullet}$ be a cohomologically constructible complex of sheaves of $\mathbb Q$-vector spaces
on $|X|$, with Verdier dual sheaf $\mathbf T\b$.  Suppose also that $H^r(|X|;\mathbf A\b) = 0$ for all $r$.
Suppose $K, L \subset X$ are supplementary (cf. Lemma \ref{lem-P1}) subcomplexes.
Then for each $i \in \ZZ$, Verdier duality induces an isomorphism
\[  H^{-i-1}(|L|;\t) \cong \Hom(H^{i}(|K|;\s),
\mathbb Q).
\]

\end{prop}

\proof  This may be proven using the method of \cite{Bjorner3}.  
The following proof uses the technology of the derived category of
sheaves on a stratified space, see for example, \cite{Iverson},
\cite{GelfandManin}, \cite{IH2}, \cite{Dimca}.  
Let $j:|K| \to |X|$ and $\iota:U = |X| - |K| \to |X|$
denote the inclusions of these complementary closed and open sets respectively.  
There are dual exact triangles in the derived category of sheaves on $|X|$,
\begin{diagram}[size=2em]
R\iota_{!}\iota^* \s & \rTo &&\s &
\text{  }& Rj_*j^{!} \t & \rTo && \t \\
& \luTo(1,2)_{\scriptstyle{[1]}} & \ldTo(2,2)& &&& \luTo(1,2)_{\scriptstyle{[1]}} &\ldTo(2,2) \\
& Rj_*j^*\s &&&&& R\iota_*\iota^*\t
\end{diagram}
where $i_{!}$ denotes the push-forward with proper supports, which amounts in
this case, to extension by zero. Similarly, $j^!\t$ is the sheaf of sections of $\t$ with support in $|K|$.

These triangles give rise to dual exact cohomology sequences.  In the following diagram, 
each cohomology group on the first line of (\ref{eqn-diag})
is dual to the group immediately below it on the second line: the second sequence is obtained from the first by
$\Hom(\cdot, \mathbb Q).$
\begin{equation}\label{eqn-diag}
\begin{diagram}
 H^i(|X|;\s) &\rTo& H^i(|K|;\s) &\rTo^{[1]}& H^{i+1}(|X|,|K|;\s) &\rTo&
H^{i+1}(|X|;\s) \\
H^{-i}(|X|;\t) &\lTo& H^{-i}_{|K|}(\t) &\lTo^{[1]}&
H^{-i-1}(U;\t) &\lTo& H^{-i-1}(|X|;\t)
\end{diagram}\end{equation}
Since $\s$ is acyclic, the sheaf $\t$ is also, and the connecting homomorphisms
denoted $[1]$ are isomorphisms, which implies that $H^i(|K|;\s)$ is dual to
$H^{-i-1}(U;\t).$

Finally we claim that the retraction $U = |X|-|K| \to |L|$ induces an
isomorphism on sheaf cohomology.  Let $h:|X| \to
[0,1]$ be the simplicial projection of Lemma \ref{lem-P1}.  Then
$h(U) =[0,1)$ is a half-open interval, and the complex $Rh_*\t$ is
(cohomologically) constructible on the closed interval $[0,1].$  Moreover,
we have natural isomorphisms, $H^*(U;\t) \cong H^*([0,1), Rh_*\t)$ and
$H^*(|K|;\t) \cong H^*(\left\{ 0 \right\}; Rh_*\t).$
So it suffices to show that for any cohomologically constructible complex of
sheaves $\mathbf A\b$ ($=Rh_*(\t)$ in this case) on $[0,1)$ the inclusion $\left\{0\right\} \to
[0,1)$ induces an isomorphism on cohomology, 
\begin{equation}\label{eqn-spectral-sequence}
H^*([0,1), \mathbf A\b)\cong H^*(\{0\}, \mathbf A\b| \{0\}).\end{equation} 
This is a standard fact (see, e.g. \cite{Borel} Lemma 3.8, or \cite{IH2} \S 1.4)
whose proof we include for completeness, as follows.

Modifying the preceding notation, let $j:\{0\} \to [0,1]$ (resp.~$\iota: (0,1] \to [0,1]$) denote the
closed (resp.~open) inclusion.
Using the first of the above exact triangles, equation (\ref{eqn-spectral-sequence})
is equivalent to the statement that $H^*(R\iota_!\iota^*\mathbf A^{\bullet})
= 0$ for all $i,$.
There is a spectral sequence for this cohomology group, whose
$E^2$ page is $H^i([0,1); R\iota_{!}\mathbf H^j(\iota^*\mathbf A^{\bullet}).$
Since the cohomology sheaf $\mathbf{H}^j(\iota^*\mathbf A^{\bullet})$ is constant on
the open interval $(0,1)$, it suffices to show that
$H^i([0,1); \iota_!\underline{\mathbb Q}_{(0,1)}) = 0$ for all $i,$ which 
is equivalent (again using the first triangle with $\s = 
\underline{\mathbb Q}_{[0,1)}$) to
the obvious statement that $H^i([0,1);\mathbb Q) \to
H^i(\left\{0\right\};\mathbb Q)$ is an isomorphism.  \qed

\quash{
\subsection{Remark}  If assumption (A) above is replaced by \begin{enumerate}
\item[($\text{A}^{\prime}$)] The cohomology group $H^i(|\Delta|;\s) = \begin{cases}
\mathbb Q &\text{for } i = 0,m\\ 0 &\text{otherwise} \end{cases}$
\end{enumerate}
then the same is true for the complex of sheaves $\t$, and the conclusion
of the theorem remains valid, provided the cohomology
$H^i(|K|;\s)$ is replaced by the reduced cohomology,
\[ \widetilde{H}^i(|K|,\s) = \begin{cases}
\coker{H^0(|\Delta|,\s) \to H^0(|K|,\s)} &\text{if } i=0\\
H^i(|K|,\s)&\text{otherwise} \end{cases}\]
and similarly for $\widetilde{H}^i(U,\t).$
}

\subsection{The sheaves $\mathbf{A}_J$}\label{subsec-AJ}  Let $\Delta $ denote the $n-1$-simplex.
  If $\sigma$ is a simplex denote its interior by $\sigma^o$
and let $i:\Delta^o \to |\Delta|$ be the inclusion.  We also denote by
\[ \QQ_! = i_!(\QQ_{\Delta^o})\]
the constant sheaf on the interior $\Delta^o$ extended by zero to all of $|\Delta|$.

The dualizing complex $\omega_{\Delta}$ on the simplex $\Delta = \Delta^{n-1}$ is a complex of sheaves whose
stalk cohomology at any point $x \in |\Delta|$ is 
\[H^{-r}_x(\omega_{\Delta}) = H_r(U, U-x) = \begin{cases} \QQ & \text{ if } x \in \Delta - \partial \Delta \text{ and }
r = n-1\\ 0 & \text{ otherwise}
 \end{cases} \]
where $U\subset \delta$ is a sufficiently small regular neighborhood of the point $x$.  Therefore in the derived
category of sheaves, the (shifted) dualizing complex is represented by the single sheaf,
\[\omega_{\Delta}[-n+1] = \QQ_!\]
where $[t]$ denotes the shift $(\mathbf B\b[t])^s = \mathbf B^{s+t}$.

If $J \subset [n]$ let $\alpha: |\Delta| - |\partial_J\Delta| \to |\Delta|$ denote the inclusion and
\[\mathbf A_J = \ker \left(\QQ_{\Delta} \to (i_J)_*\QQ_{\partial_J\Delta}\right) \cong \alpha_!(\QQ_{\Delta - \partial_J\Delta})\]
denote the constant sheaf on $|\Delta|- |\partial_J\Delta|$ extended by zero to all of $|\Delta|$.

\begin{lem}\label{lem-Verdier-dual} If $I,J \subset [n]$ are complementary subsets 
then the sheaves  $\mathbf A_J$ and
$\mathbf A_I[n-1]$ are Verdier dual:  
\[\mathbf{A}_I \cong {\rm {\bf RHom}}(A_J, \omega_{\Delta})[-n+1] = {\rm {\bf RHom}}
(\mathbf{A}_J, \QQ_!).\]
\end{lem}

\quash{
\subsection{}
If we write $A_J = \alpha_!\alpha^*(\QQ_{\Delta})$ then its Verdier dual is 
\[\alpha_*\alpha^!(\omega_{\Delta})= \alpha_*\alpha^*(\omega_{\Delta})\]
since $\alpha$ is the inclusion of an open subset.  Since $\omega_{\Delta}$ is zero on $\partial \Delta$,
this sheaf $\alpha_*\alpha^*(\omega_{\Delta})$ is zero at all points in the boundary except possibly those
points in $\partial \Delta$.  It is not difficult to see that the stalk cohomology is then
\[ H^{-r}_x(\alpha_*\alpha^* \omega_{\Delta}) = \begin{cases}
\QQ & \text{ if } x \in \text{int}(\partial_I\Delta)\\
0 & \text{ otherwise}
\end{cases}\]
so it coincides with $\mathbf A_I[n-1]$.  Here, $\text{int}(\partial_J\Delta)$ denotes the (relative) interior points of the union
of faces $\partial_J\Delta$, so that $\partial \Delta = \partial_I \Delta \coprod \text{int}(\partial_J\Delta)$.
}

\subsection{}  A quick but rather opaque proof may be given using technology of the derived category.
Rather, we give a more explicit proof within the category of simplicial sheaves on the simplicial complex $\Delta'$ (the
barycentric subdivision of $\Delta$).  Because the higher cohomology of these
sheaves vanish, the sheaf
$ {\rm {\bf RHom}}(A_J, \QQ_!)$ 
is quasi-isomorphic to the simplicial sheaf ${\mathbf {Hom}}(A_J,\QQ_!)$ which may be calculated directly. 
As a simplicial sheaf, 
\[ \QQ_!(\sigma) = \begin{cases} 0 & \text{ if } \sigma \subset \partial\Delta \\
\QQ & \text{ otherwise}\end{cases} \ \text{ and }\ 
 A_J(\sigma) = \begin{cases} 0 &\text{ if } \sigma \subset \partial_J\Delta \\
\QQ & \text{ otherwise}\end{cases}\]
for each simplex $\sigma \in \Delta'$.  Recall the simplicial sheaf ${\mathbf{Hom}}(\mathbf A, \mathbf B)$ assigns to
any simplex $\sigma$ the group $\Hom(\mathbf{A}|St^o(\sigma), \mathbf{B}|St^o(\sigma))$ 
where $St^o(\sigma)$ denotes the open star of $\sigma$.

Let us divide $|\partial\Delta|$ into three disjoint subsets,
\[ |\partial \Delta| = \left(|\partial_I\Delta| \cap |\partial_J\Delta|\right) \coprod \partial_I\Delta^o \coprod \partial_J\Delta^o\]
where $\partial_I\Delta^o = |\partial_I\Delta| - \left(|\partial_I\Delta| \cap |\partial_J\Delta|\right)$.
 If $|\sigma| \subset |\partial \Delta|$ we need to show that
\begin{equation}\label{eqn-evaluate}
\ {\mathbf{Hom}}(A_J, \QQ_!)(\sigma) =\begin{cases} 0 &\text{ if } |\sigma| \subset |\partial_I\Delta|\\
\QQ &\text{ otherwise}\end{cases} = 
 \begin{cases}
0 & \text{ if } \sigma^o \subset \partial_I\Delta \cap \partial_J\Delta\\
0 &\text{ if } \sigma^o \subset \partial_I\Delta^o\\
\QQ &\text{ if } \sigma^o \subset \partial_J\Delta ^o
\end{cases}.\end{equation}
If $J = \phi$ or $I= \phi$ these statements are easily checked so suppose that $|\partial_I\Delta| \cap |\partial_J\Delta| \ne \phi$.
Let us check the first condition in (\ref{eqn-evaluate}).
If $|\sigma| \subset |\partial_I\Delta| \cap |\partial_J\Delta|$ then there is a simplex $\tau>\sigma$
with $\tau^o\subset |\partial_I\Delta^o|$.  Let
$\widetilde{\tau}, \widetilde{\sigma}$ denote the join with the interior vertex $\ul{\mathbf 1} \in \Delta^o$
so that

\begin{diagram}[height=1.5em]
\widetilde{\tau} & > & \widetilde{\sigma}  \\
\vee && \vee  \\
\tau_1 & > & \sigma 
\end{diagram}
Then an element $\phi\in{\rm \bf Hom}(A_J, \QQ_!)(\sigma)$ is determined by a commutative mapping
\[A_J(*) \overset{\phi}{\longrightarrow} \omega_{\Delta}(*),\] that is,
\begin{diagram}[size=2em]
\QQ & \cong & \QQ   && \QQ & \cong& \QQ  \\
\uTo^{\cong} && \uTo                   &\quad \rTo^{\phi}\quad & \uTo && \uTo  \\
\QQ & \lTo & 0  &&  0 & \lTo & 0 
\end{diagram}
So $\phi = 0$.  
The other two conditions in (\ref{eqn-evaluate}) are similar but easier.  \qed

\end{document}